\documentclass[10pt]{amsart}

\usepackage{amsmath,amssymb,amscd,mathrsfs}

\theoremstyle{definition}
\newtheorem{thm}{Theorem}[section]
\newtheorem{lem}[thm]{Lemma}
\newtheorem{prp}[thm]{Proposition}
\newtheorem{dfn}[thm]{Definition}
\newtheorem{cor}[thm]{Corollary}

\newtheorem{cnj}[thm]{Conjecture}

\newtheorem{rmk}[thm]{Remark}
\newtheorem{ntn}[thm]{Notation}

\newenvironment{pff}{{\em Proof:}}{\QED}

\newcommand{\beq}{\begin{equation}}
\newcommand{\eeq}{\end{equation}}
\newcommand{\beqr}{\begin{eqnarray*}}
\newcommand{\eeqr}{\end{eqnarray*}}

\newcommand{\bit}{\begin{itemize}}
\newcommand{\eit}{\end{itemize}}

\newcommand{\af}{\alpha}
\newcommand{\bt}{\beta}

\newcommand{\dt}{\delta}

\newcommand{\ld}{\lambda}
\newcommand{\sm}{\sigma}

\newcommand{\ph}{\varphi}

\newcommand{\Ld}{\Lambda}
\newcommand{\Sm}{\Sigma}

\newcommand{\Om}{\Omega}

\newcommand{\uld}{\underline{\ld}}
\newcommand{\umu}{\underline{\mu}}
\newcommand{\uaf}{\underline{\af}}
\newcommand{\ubt}{\underline{\bt}}

\newcommand{\ui}{\underline{i}}
\newcommand{\uj}{\underline{j}}

\newcommand{\hG}{\widehat{G}}

\newcommand{\Z}{{\mathbb{Z}}}

\newcommand{\C}{{\mathbb{C}}}

\newcommand{\Sl}{{\mathfrak{sl}}}

\newcommand{\g}{{\mathfrak{g}}}
\newcommand{\h}{{\mathfrak{h}}}
\newcommand{\T}{{\mathfrak{t}}}
\newcommand{\U}{{U_{\Z}}}
\newcommand{\V}{{V_{\Z}}}

\newcommand{\B}{{\mathcal{B}}}

\newcommand{\diag}{{\mathrm{diag}}}

\newcommand{\Hom}{{\mathrm{Hom}}}

\newcommand{\End}{{\mathrm{End}}}

\newcommand{\Par}{{\mathrm{Par}}}

\newcommand{\la}{{\langle}}
\newcommand{\ra}{{\rangle}}

\newcommand{\va}{{|0\ra}}

\newcommand{\andeqn}{\,\,\,\,\,\, {\mbox{and}} \,\,\,\,\,\,}
\newcommand{\QED}{\rule{0.4em}{2ex}}

\begin{document}

\begin{center}
\textbf{Elementary Divisors of the Shapovalov Form on the Basic
Representation of Kac-Moody Algebras}

\bigskip

\emph{David Hill}

\emph{University of Oregon}

\emph{Eugene, Oregon}

\emph{U.S.A.}

E-mail: \verb"dhill1@uoregon.edu"
\end{center}

\section{Introduction/Prospectus}

Let $\g$  be the simply-laced Kac-Moody algebra of type
$X_l^{(1)}$ (i.e. $X=A,D,$ or $E$) and let $V=V(\Ld_0)$ be the
basic representation of $\g$. Let $\va$ be a vacuum vector and
define the lattice $\V=\U\va$ in $V$, where $\U$ is the
$\Z$-subalgebra of the universal enveloping algebra of $\g$
generated by the divided powers
\[
\frac{e_i^n}{n!}\andeqn\frac{f_i^n}{n!},\;\;\;(i=0,\ldots,
l,\;n\geq1)
\]
in the Chevalley generators.

Let $\tau:\g\rightarrow\g$ be the antilinear Chevalley
antiautomorphism defined on the Chevalley generators by
\[
\tau(e_i)=f_i,\;\;\;\tau(f_i)=e_i.
\]
This involution extends to an involution of the universal
enveloping algebra $U(\g)$. The Shapovalov form,
$\la\cdot,\cdot\ra_S$, is the unique Hermitian form on $V$
satisfying
\[
\la\va,\va\ra_S=1 \andeqn \la xv,v'\ra_S=\la v,\tau(x)v'\ra_S
\]
for $x\in U(\g)$ and all $v,v'\in V$. The restriction of the
Shapovalov form to $\V$ gives a symmetric bilinear form
\[
\la\cdot,\cdot\ra_S:\V\times\V\rightarrow\Z.
\]
The main result of this paper is an algorithm for calculating the
invariant factors of this form. We give formulas for the invariant
factors of the Gram matrix of this form on each weight space of
$\V$, provided the powers of the primes occurring in $l+1$ are not
too large.

Indeed, we begin by proving a relationship between the Shapovalov
form on $V$ and a family of bilinear forms $\la\cdot,\cdot\ra_s$,
($s\in\Z_{\geq0}$), on the ring $\Ld_\C=\C[p_1,p_2,\ldots]$ of
symmetric functions defined by
\[
\la p_\ld,p_\mu\ra_s=\dt_{\ld\mu}s^{l(\ld)}z_\ld
\]
where $p_\ld$ is a \emph{power sum} symmetric function, and for
$\ld=(1^{m_1}2^{m_2}\cdots)$,
\[
z_{\ld}=\prod_{r\geq 1}(r^{m_r}\cdot m_r!).
\]

This form restricts to a symmetric bilinear form on
$\Ld:=\Ld_\Z=\Z[h_1,h_2,\ldots]$, where $h_n$ is the $n$th
\emph{complete homogeneous} symmetric function. There is a simple
relationship between $\V$ and $\Ld$. Indeed, it turns out that the
weight spaces of $V$ are of the form $w\Ld_0-d\dt$, where $w$ is
an element of the Weyl group of type $X$, $\dt$ is the null root,
and $d$ is a positive integer (see \cite{kc}, section 12.6). We
have

\begin{thm}Let $a_1,\ldots,a_l$ be the invariant factors of the
Cartan matrix for a simple finite dimensional Lie algebra of $ADE$
type. Let $a^{(r)}_{k1},\ldots,a^{(r)}_{kh}$ be the invariant
factors of the form $\la\cdot,\cdot\ra_{a_k}$ on the degree $r$
component of $\Ld$ (here $h=|\Par(r)|$ is the number of partitions
of $r$). Then, the invariant factors of the Shapovalov form on the
$(w\Ld_0-d\dt)$-weight space of $\V$ are
\[
\left\{\prod_{k=1}^la_{ki_k}^{(d_k)}:d_1+\cdots+d_l=d, 1\leq
i_k\leq|\Par(d_k)| \right\}.
\]
\end{thm}

This reduces the problem to calculating the invariant factors of
the form $\la\cdot,\cdot\ra_s$ on $\Ld$ for positive integers $s$.

We further reduce the problem as follows. Let $X_s$ denote the
matrix $(\la m_\ld,h_\mu\ra_s)_{\ld,\mu}$ ($m_\ld$ is a
\emph{monomial} symmetric function). The $m_\ld$ also form a basis
for $\Ld$, so we may compute the invariant factors of $X_s$. This
calculation is equivalent to calculating the Smith normal form
$S(X_s)$ of $X_s$. We have:

\begin{thm} Let $s,t\in\Z_{\geq0}$. Then, $X_{st}=X_sX_t$. In
particular, if $(s,t)=1$, then $S(X_{st})=S(X_s)S(X_t)$.
\end{thm}

Therefore, to calculate invariant factors on $\V$, it is enough to
know the invariant factors of the form $\la\cdot,\cdot\ra_{p^r}$
on $\Ld$ for every prime $p$ and $r\geq1$. We calculate these
numbers in the case when $r\leq p$.

For an integer $a$, define the number
\begin{eqnarray*}\label{formula: geometric sum}
d_p(a)=\sum_{j\geq 1}\left\lfloor\frac{a}{p^j}\right\rfloor
\end{eqnarray*}
and for a partition $\ld=(1^{m_1(\ld)}2^{m_2(\ld)}\ldots)$, define
\begin{eqnarray}\label{eqn: invariant factors}
D_r(\ld)=\prod_{(n,p)=1}\prod_{i=0}^{r-1}p^{[(r-i)m_{p^in}(\ld)+d_p(m_{p^in}(\ld))]}
\end{eqnarray}
Then,

\begin{thm}\label{thm:main}
Let $r\leq p$. Then, the elementary divisors of the form
$\la\cdot,\cdot\ra_{p^r}$ on the $d$th graded component of $\Ld$
are
\[
\{D_r(\ld) |\ld\vdash d\}.
\]
\end{thm}

We conjecture that this result holds for all $r$ based on
computational evidence, but have yet to understand the proof.

This result has a number of consequences. First, assume $\g$ is of
type $A_{l-1}^{(1)}$. Let $H_n$ be the Iwahori-Hecke algebra
associated to the symmetric group $S_n$ at a primitive $l$th root
of unity. In \cite{ar} and \cite{gr}, it was shown that there is
an isomorphism between the basic representation $\V$ of $\g$ and
the direct sum $K=\bigoplus_nK_n$ of Grothendieck groups $K_n$ of
finitely generated projective $H_n$-modules. Under this
isomorphism there is a correspondence between blocks of $K$ and
weight spaces of $\V$. The Shapovalov form on $\V$ corresponds to
the Cartan pairing
\[
([P],[Q])=\dim\Hom(P,Q)
\]
between projective modules $P$ and $Q$. As a consequence, this
paper gives the invariant factors of the blocks of the Hecke
algebra $H_n$ at an $l$th root of unity when $l=\prod_ip_i^{r_i}$
satisfies $r_i\leq p_i$ for all $i$.

In \cite{kor}, K\"{u}lshammer, Olsson and Robinson develop an
$l$-modular representation theory of symmetric groups for $l$ not
necessarily prime. They conjecture what the invariant factors of
the blocks should be. We expect that the answer given in this
paper will confirm their conjecture.

Finally, we expect that these results will extend to include
twisted affine Kac-Moody algebras.

\section{Reduction to the Heisenberg subalgebra of $\g$}

Let $\g'$ be the finite dimensional simple Lie algebra of type
$X_l$ with root system $\Phi$, simple roots
$\af_1',\ldots,\af_l'$, and root lattice
$Q=\bigoplus_{i=1}^l\Z\af_i'$. Let $\h'=\C\otimes_{\Z}Q\subset\g'$
be its Cartan subalgebra. Recall that
\[
\g=\g'\otimes\C[t,t^{-1}]\oplus\C c\oplus\C d
\]
where $c$ is the one dimensional central extension of
$\g'\otimes\C[t,t^{-1}]$ and $d$ is the scaling element. Let
$(\cdot|\cdot)$ be a standard invariant form on $\g$, normalized
so that $(\af|\af)=2$ for each root $\af\in Q$.

Notice that $\g$ contains a copy of the Heisenberg Lie algebra
$\T=\T^+\oplus\C c\oplus\T^-$, where
\[
\T^{\pm}=\bigoplus_{\pm n>0}\h'\otimes t^n.
\]
In particular, $\T$ is generated by elements $\af_i'(n)$
($i=0,\ldots,l$, $n\in\Z$), where $\af_0'$ is the longest root of
$\Phi$ and
\[
\af_i'(n)=\af_i\otimes t^n.
\]
The commutation relations for $\T$ are given by
\[
[\af_i'(n),\af_j'(m)]=n(\af_i'|\af_j')\dt_{n,-m}c.
\]

We may view the symmetric algebra $S(\T^-)$ as a $\T$-module so
that $c$ acts as 1, elements of $\T^-$ act by multiplication and
elements of $\T^+$ annihilate 1. It is $\Z$-graded by setting
\[
\deg(h\otimes t^{-n})=n
\]
for $h\in\h'$ and $n\geq1$.

Define generating series
\[
H_{\af_i'}(z)=\exp\left(\sum_{n\geq
1}\frac{\af_i'(-n)z^n}{n}\right)\andeqn
K_{\af_i'}(z)=\left(-\sum_{n\geq1}\frac{\af_i(n)'z^n}{n}\right)
\]
viewed as elements of $\End(S(\T^-))[[z^{\pm1}]]$.

Next, for $n\geq 1$ and $i=0,\ldots,l$, define
\[
y_n^{(i)}=\af_i'(-n)/n\andeqn x_n^{(i)}=\sum_{k_1+2k_2+\cdots=n}
\frac{(y_1^{(i)})^{k_1}}{k_1!}\frac{(y_2^{(i)})^{k_2}}{k_2!}
\frac{(y_3^{(i)})^{k_3}}{k_3!}\cdots.
\]
Note that
\[
H^{(i)}(t):=H_{\af_i'}(t)=\exp\left(\sum_{n\geq1}\frac{y_n^{(i)}}{n}t^n\right)=1+\sum_{n\geq1}x_n^{(i)}t^n.
\]
The $y_n^{(i)}$, $i=1,\ldots,n$, $n\geq1$, form a basis for
$\T^-$. Hence $S(\T^-)$ is the free polynomial algebra
\[
\B_{\C}:=\bigotimes_{i=1}^{l}\C[y^{(i)}_n|n\geq1].
\]
Since the transition matrix from the $x_n^{(i)}$ to the
$y^{(i)}_n/n$ is unitriangular, we obtain a $\Z$-form
\[
\B_{\Z}=\bigotimes_{i=1}^l\Z[x^{(i)}_n|n\geq 1]\subset \B_{\C}
\]
for $\B_{\C}$. Moreover, since $\af_0'$ is a $\Z$-linear
combination of simple roots, the elements $x_n^{(0)}$ belong to
the lattice $\B_\Z$. Finally, $\B_\Z$ inherits a $\Z$-grading from
the grading on $S(\T^-)$ given by
\[
\deg(y^{(i)}_n)=\deg(x^{(i)}_n)=n
\]
for $i=1,\ldots,l$.

Fix $d\geq0$. It was shown in (\cite{bk}, lemma 4.1) that the Gram
matrix of the Shapovalov form on the ($w\Ld_0-d\dt$) weight space
of $\V$ is related to the Gram matrix of the Shapovalov form on
the degree $d$ component of $\B_\Z$ in a unimodular way.

Inspiration for this paper comes from (\cite{bk}, 3.5). The
relevant fact for our purposes is that
\[
K_{\af_i'}(z)H_{\af_j'}(w)=H_{\af_j'}(w)K_{\af_i'}(z)(1-zw)^{-a_{ij}}
\]
where $a_{ij}=(\af_i'|\af_j')$ is the $(i,j)$-entry in the Cartan
matrix for $\g'$. In the next section, we will show how to extract
the entries of the Gram matrix of the Shapovalov form from certain
coefficients in the generating series
\[
\prod_{i,j=1}^l\prod_{s,t\geq 1}(1-z_s^{(i)}w_t^{(j)})^{-a_{ij}}.
\]

\section{The Shapovalov Form on the Ring of Symmetric
Functions}\label{section:symmetric functions}

Throughout the paper, we will use notation from MacDonald,
\cite{m}. In particular, for a partition $\ld$,  $m_\ld$, $h_\ld$,
and $p_\ld$ will denote the monomial, homogeneous, and power sum
symmetric functions, respectively.

Now, since $\B_\Z$ is a free polynomial algebra, we may identify
it with the $l$-fold tensor product
\[
\Ld:=\bigotimes_{i=1}^l\Ld^{(i)}
\]
via $x_n^{(i)}=h_n^{(i)}$. Here,
$\Ld^{(i)}$ is the ring of symmetric functions in the
\emph{variables} $z_n^{(i)}$ (or any other letter) and
$h_n^{(i)}(z):=h_n(z^{(i)})$ is the $n$th homogeneous symmetric
function in the variables $z_n^{(i)}$, $n\geq1$.

This induces an identification of $\mathcal{B}_{\C}$ with
$\Ld_{\C}=\bigotimes_{i=1}^l\Ld^{(i)}_{\C}$ so that
$y_n^{(i)}=p_n^{(i)}/n$.

We may index bases of $\Ld$ (resp. $\Ld_\C$) as follows:

\begin{ntn}\label{index:1} Let $I=\{1,2,\ldots,l\}$, and for a partition
$\ld=(\ld_1\geq\ld_2,\cdots\geq\ld_k>0)$, let $I(\ld)=I^{k}$. If
$\underline{i}=(i_1,\ldots,i_k)\in I(\ld)$, set
\[
h_{\ld}^{(\underline{i})}=h_{\ld_1}^{(i_1)}\cdots
h_{\ld_k}^{(i_k)},\;\;\;\;\;
p_{\ld}^{(\underline{i})}=p_{\ld_1}^{(i_1)}\cdots
p_{\ld_k}^{(i_k)},\;\;\;\;\;
m_{\ld}^{(\underline{i})}=m_{\ld^{(1)}}^{(1)}\cdots
m_{\ld^{(l)}}^{(l)}
\]
where $\ld^{(j)}$ is the partition consisting of the parts $\ld_r$
satisfy $i_r=j$. (Note the difference in the formula describing
the monomial symmetric functions. This is due to the fact that the
$m_\ld$ are not multiplicative.)

Let $\Om(\ld)=\{\underline{i}\in I(\ld)|i_j\leq i_{j+1}\mbox{
whenever }\ld_j=\ld_{j+1}\}$. Then the set $\{(\ld,\underline{i})|
\ld\in\Par(d), \underline{i}\in\Om(\ld)\}$ is in one-to-one
correspondence with the dimension of the $d$th graded component
$\Ld_d$ of $\Ld$.
\end{ntn}

It is worth noting another indexing system for bases of $\Ld$.

\begin{ntn}\label{index:2} Given $\ld\in\Par(d)$ and $\ui\in\Om(\ld)$, let $\ld^{(j)}$ be the
partition consisting of the parts $\ld_r$ satisfy $i_r=j$. Define
the \emph{multipartition}
$\uld:=(\ld^{(1)},\ld^{(2)},\ldots,\ld^{(l)})$, set
$|\uld|=\sum|\ld^{(i)}|$. Define
\[
h_{\uld}=h_{\ld^{(1)}}^{(1)}\cdots h_{\ld^{(l)}}^{(l)},\;\;\;\;\;
p_{\uld}=p_{\ld^{(1)}}^{(1)}\cdots p_{\ld^{(l)}}^{(l)},\;\;\;\;\;
m_{\uld}=m_{\ld^{(1)}}^{(1)}\cdots m_{\ld^{(l)}}^{(l)}.
\]
Then, the set of all multipartitions $\uld$ with $|\uld|=d$ is in
one-to-one correspondence with the dimension of $\Ld_d$.

When working with multipartition $\uld$, it will sometimes be
useful to refer to the associated partition $\Sm\uld$ defined by
$\Sm\uld=\sum_{i=1}^l\ld^{(i)}$. Here, given two partitions
$\sm=(i^{m_i(\sm)})$ and $\rho=(i^{m_i(\rho)})$,
$\sm+\rho=(i^{m_i(\sm)+m_i(\rho)})$.

To go back to the notation of \ref{index:1}, note that
$h_{\uld}=h_{\Sm\uld}^{(\ui)}$, where $\ui=(i_1,\ldots,i_k)$ is
obtained by writing
$\Sm\uld=(\ld^{(i_1)}_1\geq\cdots\geq\ld^{(i_k)}_k>0)$ according
to the rule that $i_j\leq i_{j+1}$ if $\ld_j=\ld_{j+1}$.
\end{ntn}

Now, we transport the Shapovalov form to $\Ld$ using the
identification above. It follows from (\cite{bk}, 5.3) that
\begin{eqnarray}\label{eqn:block diagonal}
\la p_\ld^{(\ui)},p_\mu^{(\uj)}\ra_S=\dt_{\ld\mu}a_{i_1j_1}\cdots
a_{i_kj_k}z_\ld
\end{eqnarray}
where, for $\ld=(1^{m_1}2^{m_2}\cdots)$,
\[
z_{\ld}=\prod_{r\geq 1}(r^{m_r}\cdot m_r!).
\]

We will now record some facts about transition matrices from
\cite{m}. Let $L:=M(p,m)$ denote the transition matrix from the
power sums to the monomial symmetric functions. Using Notation
\ref{index:1}, $L=(L_{\uld\umu})$ where
\[
p_{\uld}=\sum_{\umu}L_{\uld\umu}m_{\umu}.
\]
The transition matrix $A=M(h,p)$ is related to $L$ by
\[
A=L^tz^{-1}
\]
where $z=\diag(z_{\Sm\uld})$ is diagonal and $\Sm\uld$ is defined
in \ref{index:2}.

We also have the unimodular transition matrix $N=M(h,m)$. Using
standard properties of transition matrices (\cite{m}, Ch.I, (6.3))
we deduce that \beqr
A   &=& M(h,p)\\
    &=& M(h,m)M(m,p)\\
    &=& NL^{-1}.
\eeqr

Let $C=(\la h_{\uld},h_{\umu}\ra_S)_{\uld\umu}$ and let $B$ be the
matrix
\[
B=(\la p_{\uld},z_{\Sm\umu}^{-1}p_{\umu}\ra_S)_{\uld\umu}.
\]
This matrix is block diagonal in the sense that $B_{\uld\umu}=0$
unless $\Sm\uld=\Sm\umu$ (see (\ref{eqn:block diagonal}) above).
Then, by the relevant definitions,
\[
A^{-1}C(A^{-1})^t=Bz
\]
whence
\[
C=ABzA^t=NL^{-1}Bzz^{-1}L
\]
so that
\[
N^{-1}C=L^{-1}BL.
\]

Notice that $N$ is a unimodular matrix, so the elementary divisors
of $N^{-1}C$ are the same as those for $C$. Notice too that

\beqr
(N^{-1}C)_{\uaf,\ubt}&=&\sum_{\uld}N^{-1}_{\uaf\uld}C_{\uld\ubt}\\
        &=&\sum_{\uld}N^{-1}_{\uaf\uld}\langle h_{\uld},h_{\ubt}\rangle_S\\
        &=&\langle\sum_{\uld}N^{-1}_{\uaf\uld}h_{\uld},h_{\ubt}\rangle_S\\
        &=&\langle m_{\uaf},h_{\ubt}\rangle_S.
\eeqr

We are now in a position to calculate a generating series for the
entries of $N^{-1}C=(\langle
m_{\uld},h_{\umu}\rangle_S)_{\uld,\umu}$. In fact, we prove a more
general proposition that will be useful in the next section. For
this Proposition, we use notation \ref{index:1}. In this notation,
$L=(L_{\ld\mu}^{\ui,\uj})$ where
\[
p_\ld^{(\ui)}=\sum_{\mu,\uj\in\Om(\mu)}L_{\ld\mu}^{\ui,\uj}m_\mu^{(\uj)}.
\]

\begin{prp}\label{prp:generating series}
Let $\tilde{P}=(p_{ij})$ be an $l\times l$ matrix and define
\[
P=\bigoplus_{\ld\in\Par}S^{m_1(\ld)}(\tilde{P})\otimes
S^{m_2(\ld)}(\tilde{P})\otimes\cdots.
\]
Then $(L^{-1}PL)_{\af\bt}^{\ui,\uj}$ is the coefficient of
$h_\af^{(\ui)}(x)m_{\bt}^{(\uj)}(y)$ in the product
\[
\prod_{i,j=1}^l\prod_{s,t}\left(1-x_s^{(i)}y_t^{(j)}\right)^{-p_{ij}}.
\]
\end{prp}

\begin{pff}
First, we describe what we mean by $S^m(\tilde{P})$. To this end,
label the rows (resp. columns) of $\tilde{P}$ with integers
$1,2,\ldots,l$ from left to right (resp. top to bottom). Then, the
rows (resp. columns) of $S^m(\tilde{P})$ are labelled by
$m$-tuples $(i_1\leq i_2\leq\cdots\leq i_m)=:\ui$, and the
$(\ui,\uj)$ entry of $S^m(\tilde{P})$ is
\[
p_{i_1j_1}p_{i_2j_2}\cdots p_{i_mj_m}.
\]
Notice here that the rows and columns of the $\ld$ block of $P$,
\[
S^{m_1(\ld)}(\tilde{P})\otimes
S^{m_2(\ld)}(\tilde{P})\otimes\cdots,
\]
are naturally labelled by $\Om(\ld)$.

We deduce the following from (\cite{s}, Proposition 7.7.4)
\begin{eqnarray*}
\prod_{i,j=1}^l\prod_{s,t\geq1}(1-x_s^{(i)}y_t^{(j)})^{-p_{ij}}
    &=&\prod_{i,j=1}^l\exp\left(\sum_{n\geq 0}\frac{1}{n}p_n^{(i)}(x)p_n^{(j)}(y)\right)^{p_{ij}}\\
    &=&\prod_{i,j=1}^l\prod_{n\geq0}\exp\left(\frac{p_{ij}}{n}
    p_n^{(i)}(x)p_n^{(j)}(y)\right)\\
    &=&\sum_{\ld,\underline{i},\underline{j}}P_{\ld\ld}^{\underline{i},\underline{j}}z_{\ld}^{-1}
p_{\ld}^{(\underline{i})}(x)p_{\ld}^{(\underline{j})}(y).
\end{eqnarray*}
Moreover, we have that
\[
(L^{-1}PL)_{\af\bt}^{\underline{k},\underline{l}}=
\sum_{\ld,\underline{i},\underline{j}}(L^{(-1)})_{\af\ld}^{\underline{k},\underline{i}}
P_{\ld\ld}^{\underline{i},\underline{j}}L_{\ld\bt}^{\underline{j},\underline{l}}.
\]
Thus,
\beqr
\prod_{i,j=1}^l\prod_{s,t\geq1}(1-x_s^{(i)}y_t^{(j)})^{-p_{ij}}
&=&\sum_{\ld,\underline{i},\underline{j}}P_{\ld\ld}^{\underline{i},\underline{j}}z_{\ld}^{-1}
p_{\ld}^{\underline{i}}(x)p_{\ld}^{\underline{j}}(y)\\
&=&\sum_{\ld,\underline{i},\underline{j}}P_{\ld\ld}^{(\underline{i})(\underline{j})}z_{\ld}^{-1}
\left(\sum_{\af,\underline{k}}(L^{(-1)})_{\af\ld}^{\underline{k},\underline{i}}z_{\ld}h_{\af}^{(\underline{k})}(x)\right)\\
&&\;\;\;\;\;\times\left(\sum_{\bt,\underline{l}}L_{\ld\bt}^{\underline{j},\underline{l}}m_{\bt}^{(\underline{l})}(y)\right)\\
&=&\sum_{\af,\underline{k},\bt,\underline{l}}\left(\sum_{\ld,\underline{i},\underline{j}}
(L^{(-1)})_{\af\ld}^{\underline{k},\underline{i}}P_{\ld\ld}^{\underline{i},\underline{j}}
L_{\ld\bt}^{\underline{j}\underline{l}}\right)h_{\af}^{(\underline{k})}(x)m_{\bt}^{(\underline{l})}(y).
\eeqr
\end{pff}

\begin{cor}\label{cor:Shapovalov generating series} The coefficient of $h_{\uaf}(x)m_{\ubt}(y)$
in the generating series
\[
\prod_{i,j=1}^l\prod_{s,t}\left(1-x_s^{(i)}y_t^{(j)}\right)^{-a_{ij}}
\]
is $\langle m_{\uaf},h_{\ubt}\rangle_S$.
\end{cor}

\begin{pff} This follows from the fact that
\[
B=\bigoplus_{\ld=(1^{m_1}2^{m_2}\cdots)}S^{m_1}(X)\otimes
S^{m_2}(X)\otimes\cdots
\]
where $X=(a_{ij})_{i,j=1}^l$ is the Cartan matrix of the simple
Lie algebra of type $X_l$ ($X=ADE$). See (\cite{bk}, 5.3) for
details.
\end{pff}

\section{Some Reductions}

\subsection{Block Diagonalization of the Gram
Matrix}\label{subsection:block diagonal} Let $X=(a_{ij})$ be the
Cartan matrix for the simple finite dimensional Lie algebra of
type $X_l$ ($X=A,D,E$). We will show that the elementary divisors
of the Shapovalov form depend only on the elementary divisors of
$X$. To this end, let $\tilde{Q}=(q_{ij})$ and
$\tilde{T}=(t_{ij})$ be unimodular matrices such that
$\tilde{Q}X\tilde{T}=\diag(a_1,a_2,\ldots,a_l)$. Construct
matrices $Q$ and $T$ by the formulae:
\[
Q=\bigoplus_{\ld\in\Par}S^{m_1(\ld)}(\tilde{Q})\otimes
S^{m_2(\ld)}(\tilde{Q})\otimes\cdots\andeqn
T=\bigoplus_{\ld\in\Par}S^{m_1(\ld)}(\tilde{T})\otimes
S^{m_2(\ld)}(\tilde{T})\otimes\cdots
\]
where $\ld=(1^{m_1}2^{m_2}\cdots)$. Then, there exist bases
$(q_{\ld})$ and $(t_{\ld})$ for $\Ld_\C$ defined by $Q=M(q,p)$ and
$T=M(t,z^{-1}p)^t$. Note that
\[
QBT=\bigoplus_{\ld}S^{m_1}(\tilde{Q}X\tilde{T})\otimes
S^{m_2}(\tilde{Q}X\tilde{T})\otimes\cdots
\]
so
\[
\langle
q_{\underline{\ld}},t_{\underline{\mu}}\rangle_{S}=\dt_{\underline{\ld}\underline{\mu}}a_1^{l(\ld^{(1)})}\cdots
a_l^{l(\ld^{(l)})}
\]
where $\underline{\ld}=(\ld^{(1)},\ldots,\ld^{(l)})$ is a
multipartition.

Define bases $(a_{\underline{\ld}})$ and $(b_{\underline{\ld}})$
for $\Ld_{\C}$ by the formulae:
\[
Y:=M(a,h)^t=L^{-1}QL\andeqn Z:=M(b,m)=L^{-1}TL.
\]
From Proposition \ref{prp:generating series} and the fact that
$q_{ij}$ and $t_{ij}$ are integers it follows that $Y$ and $Z$ are
integral matrices. As their determinant is clearly 1, they are
unimodular. Hence $(a_{\uld})$ and $(b_{\uld})$ are bases for
$\Ld$.

\begin{prp}\label{prp:basis change}
\[
\prod_{i,j}\prod_{k=1}^l(1-x_i^{(k)}y_j^{(k)})^{-a_k}=
\sum_{\underline{\af},\underline{\bt}}\langle
m_{\underline{\af}},h_{\underline{\bt}}\rangle_S
a_{\underline{\af}}b_{\underline{\bt}}.
\]
\end{prp}

\begin{pff}
Let
\[
R:=M(z^{-1}p,a)^t=L^{-1}Q^{-1}\andeqn S:=M(p,b)=T^{-1}L.
\]
Then,
\[
N^{-1}C=L^{-1}BL=(L^{-1}Q^{-1})(QBT)(T^{-1}L)=R(QBT)S.
\]
Therefore,
\begin{eqnarray*}
\langle
m_{\underline{\af}},h_{\underline{\bt}}\rangle_S&=&(N^{-1}C)_{\underline{\af}\underline{\bt}}\\
    &=&\sum_{\underline{\ld}}a_1^{l(\ld^{(1)})}\cdots
    a_l^{l(\ld^{(l)})}R_{\underline{\af}\underline{\ld}}S_{\underline{\ld}\underline{\bt}}.
\end{eqnarray*}
Now, by a calculation similar to the one given in Proposition
\ref{prp:generating series}, one has
\begin{eqnarray*}
\prod_{i,j}\prod_{k=1}^l(1-x_i^{(k)}y_j^{(k)})^{-a_k}
    &=&\prod_{k=1}^l\prod_{n\geq0}\exp\left(\frac{a_k}{n}p_n^{(k)}(x)p_n^{(k)}(y)\right)\\
    &=& \sum_{\underline{\ld}}a_1^{l(\ld^{(1)})}\cdots
a_l^{l(\ld^{(l)})}
z_{\underline{\ld}}^{-1}p_{\underline{\ld}}(x)p_{\underline{\ld}}(y)\\
    &=&\sum_{\underline{\ld}}a_1^{l(\ld^{(1)})}\cdots a_l^{l(\ld^{(l)})}
        \left(\sum_{\underline{\af}}R_{\underline{\af}\underline{\ld}}a_{\underline{\af}}(x)\right)\\
        &&\hspace{.75in}\times\left(\sum_{\underline{\bt}}S_{\underline{\ld}\underline{\bt}}b_{\underline{\bt}}(y)\right)\\
    &=&\sum_{\underline{\af},\underline{\bt}}\left(\sum_{\underline{\ld}}a_1^{l(\ld^{(1)})}\cdots a_l^{l(\ld^{(l)})}
        R_{\underline{\af}\underline{\ld}}S_{\underline{\ld}\underline{\bt}}\right)
        a_{\underline{\af}}(x)b_{\underline{\bt}}(y)\\
    &=&\sum_{\underline{\af},\underline{\bt}}\langle
        m_{\underline{\af}},h_{\underline{\bt}}\rangle_S
        a_{\underline{\af}}(x)b_{\underline{\bt}}(y).
\end{eqnarray*}
\end{pff}

Now, if we define bases $(c_{\underline{\ld}})$ and
$(d_{\underline{\ld}})$ for $\Ld$ by the formulae $M(c,m)=Y$ and
$M(d,h)^t=Z$, we obtain
\[
\sum_{\underline{\af},\underline{\bt}}\langle
        m_{\underline{\af}},h_{\underline{\bt}}\rangle_S
        a_{\underline{\af}}(x)b_{\underline{\bt}}(y)=
\sum_{\underline{\af},\underline{\bt}}\langle
        c_{\underline{\af}},d_{\underline{\bt}}\rangle_S
        h_{\underline{\af}}(x)m_{\underline{\bt}}(y).
\]
Since the $c_{\underline{\ld}}$ and $d_{\underline{\ld}}$ are
obtained from $h_{\underline{\ld}}$ and $m_{\underline{\ld}}$ by
unimodular change, the matrix
\[
YN^{-1}CZ=(\langle
c_{\underline{\af}},d_{\underline{\bt}}\rangle_S)_{\underline{\af},\underline{\bt}}
\]
has the same elementary divisors as $N^{-1}C$.

\subsection{From the MacDonald Pairing to the Shapovalov Form} From now on,
let $\Ld$ (resp. $\Ld_\C$) be the usual ring of symmetric
functions over $\Z$ (resp. $\C$). Given any real number $s>0$,
consider the scalar product $\langle\cdot,\cdot\rangle_s$ on
$\Ld_{\C}$ defined on the power sum symmetric functions by the
formula
\[
\langle p_{\ld},p_{\mu}\rangle_s=\dt_{\ld\mu}s^{l(\ld)}z_{\ld}
\]
as in (\cite{m}, Ch.VI, section 10). We call this form the
$s$-form on $\Ld_\C$.

Since the following term occurs frequently, we make the
abbreviation:
\[
\Pi(x,y):=\prod_{i,j}(1-x_iy_j)^{-1}.
\]

\begin{prp}\label{prp:s-form} If $(u_{\ld})$ and $(v_{\ld})$ are two bases for
$\Ld$, then
\[
\Pi(x,y)^{s}
    =\sum_{\ld,\mu}\langle u_{\ld}^*,v_{\mu}^*\rangle_s
    u_{\ld}(x)v_{\mu}(y)
\]
where $(u_{\ld}^*)$ (resp. $(v_{\ld}^*)$) is the dual basis to
$(u_\ld)$ (resp. $(v_\ld)$) with respect to the form
$\la\cdot,\cdot\ra_1$ on $\Ld$.
\end{prp}

\begin{pff}
Write $M(u,h)=(a_{\ld\mu})_{\ld,\mu}$ and
$M(v,m)=(b_{\ld\mu})_{\ld,\mu}$. Then, $M(u^*,m)=(M(u,h)^{-1})^t$
and $M(v^*,h)=(M(v,m)^{-1})^t$ (see \cite{m}, Ch.I, 6.3). Thus,

\beqr \sum_{\ld,\mu}\langle u^*_{\ld},v^*_{\mu}\rangle_s
u_{\ld}(x)v_{\mu}(y)
    &=&\sum_{\ld,\mu}\left\langle\sum_{\af}a_{\af\ld}^{(-1)}m_{\af},\sum_{\bt}b_{\bt\mu}^{(-1)}h_{\bt}\right\rangle_s\\
    &&\hspace{.75in}\times\left(\sum_{\sm}a_{\ld\sm}h_{\sm}(x)\right)\left(\sum_{\rho}b_{\mu\rho}m_{\rho}(y)\right)\\
    &=&\sum_{\af,\bt,\sm,\rho}\left(\sum_{\ld}a_{\af\ld}^{(-1)}a_{\ld\sm}\right)\left(\sum_{\mu}b_{\bt\mu}^{(-1)}b_{\mu\rho}
    \right)\\
    &&\hspace{.75in}\times\langle m_{\af},h_{\bt}\rangle_sh_{\sm}(x)m_{\rho}(y)\\
    &=&\sum_{\af,\bt}\langle
    m_{\af},h_{\bt}\rangle_sh_{\af}(x)m_{\bt}(y)\\
    &=&\Pi(x,y)^{s}.
\eeqr
The last equality follows from Proposition
\ref{prp:generating series} by taking the matrix $\tilde{P}$ to be
the $1\times1$ matrix $(s)$.
\end{pff}

\begin{thm}\label{cor:invariant factors} Let $a_1,\ldots,a_l$ be the invariant factors of the
Cartan matrix for a simple finite dimensional Lie algebra of $ADE$
type. Let $a^{(r)}_{k1},\ldots,a^{(r)}_{kh}$ be the invariant
factors of the form $\la\cdot,\cdot\ra_{a_k}$ on the degree $r$
component of $\Ld$ (here $h=|\Par(r)|$ is the number of partitions
of $r$). Then, the invariant factors of the Shapovalov form on the
$(w\Ld_0-d\dt)$-weight space of $\V$ are
\[
\left\{\prod_{k=1}^la_{ki_k}^{(d_k)}:d_1+\cdots+d_l=d, 1\leq
i_k\leq|\Par(d_k)| \right\}.
\]
\end{thm}

\begin{pff}
The coefficient of $h_{\uld}(x)m_{\umu}(y)$ in the product
\[
\prod_{i,j}\prod_{k=1}^l(1-x_i^{(k)}y_j^{(k)})^{-a_k}
\]
is $\langle c_{\uld},d_{\umu}\rangle_S$, where $c_{\uld}$ and
$d_{\uld}$ are dual to $a_{\uld}$ and $b_{\umu}$ with respect to
the 1-form on $\Ld$ (see Proposition \ref{prp:basis change} and
subsequent remarks). Since $c_{\uld}$ and $d_{\umu}$ are obtained
from $m_{\uld}$ and $h_{\umu}$ by unimodular change, the matrix
$(\langle c_{\uld},d_{\umu}\rangle_S)$ has the same invariant
factors as $C$.

On the other hand, by Proposition \ref{prp:s-form}, the
coefficient of $h_{\uld}(x)m_{\umu}(y)$ is
\[
\prod_{k=1}^l\la m_{\ld^{(k)}}, h_{\mu^{(k)}}\ra_{a_k}.
\]
In particular,
\[
(\langle
c_{\uld},d_{\umu}\rangle_S)_{\uld,\umu}=\left(\prod_{k=1}^l\la
m_{\ld^{(k)}}, h_{\mu^{(k)}}\ra_{a_k}\right)_{\uld,\umu}.
\]

It is just left to observe that the invariant factors of the
matrix on the right hand side of the equality above are
\[
\left\{\prod_{k=1}^la_{ki_k}^{(d_k)}:d_1+\cdots+d_l=d, 1\leq
i_k\leq|\Par(d_k)| \right\}.
\]
\end{pff}

\subsection{Splitting the Gram Matrix across Primes}
For each positive integer $s$, let $X_s=(\la
m_{\uld},h_{\umu}\ra_s)_{\uld,\umu}$ denote the Gram matrix of the
$s$-form on $\Ld$, and
\[
B_s=\diag\{s^{l(\ld)}\}.
\]
Then, it follows from the calculation done in Proposition
\ref{prp:generating series} (with the matrix $(a_{ij})$ taken to
be the $1\times1$ matrix $(s)$) that
\[
X_s=L^{-1}B_sL.
\]
Notice that
\begin{eqnarray*}
X_{st}&=&L^{-1}B_{st}L\\
    &=&L^{-1}B_sB_tL\\
    &=&(L^{-1}B_sL)(L^{-1}B_tL)\\
    &=&X_sX_t.
\end{eqnarray*}
In particular, if $s$ and $t$ are positive integers satisfying
$(s,t)=1$, then
\[
(\det X_s, \det X_t)=1.
\]
Indeed, this follows since $\det(X_s)$ is a power of $s$.  By
(\cite{n}, Theorem II.15),
\[
S(X_{st})=S(X_s)S(X_t),
\]
where $S(X)$ is the Smith normal form
of a matrix $X$. This reduces the problem to calculating the
elementary divisors of $X_{p^r}=(X_p)^r$.

\section{The Invariant Factors}

\subsection{Higher Homogeneous Functions}
In the previous sections, we have reduced the problem of
calculating the invariant factors of the Shapovalov form to
calculating the invariant factors of the Gram matrix of the
$p^r$-form on the ring of symmetric functions. Thanks to
Proposition \ref{prp:s-form}, we have identified the entries in
this matrix as coefficients of the generating series
$\Pi(x,y)^{p^r}$. We have the following:

\begin{lem}\label{lem:s-form generating series}
\[
\Pi(x,y)^{s}=\sum_{\mu}\prod_{i=1}^{l(\mu)}\left(\sum_{\ld\vdash\mu_i}
    \left.{s}\choose{l(\ld)}\right. \left.{l(\ld)}\choose{m_1(\ld),m_2(\ld),\ldots}\right.
    h_{\ld}(x)\right)m_{\mu}(y).
\]
\end{lem}

\begin{pff}
This is a straightforward calculation.
\begin{eqnarray*}
\prod_{i}(1-x_iy_j)^{-s}&=&\left(1+\sum_{n\geq 1}h_n(x)y_j^n\right)^s\\
        &=&\sum_{k\geq0}\left.{s}\choose{k}\right.\left(\sum_{n\geq
        1}h_n(x)y_j^n\right)^k\\
        &=&\sum_{n\geq0}\sum_{k\geq0}\left.{s}\choose{k}\right.\left(\sum_{\af}h_{\af}(x)\right)y_j^n
\end{eqnarray*}
where the last sum is over all $k$-tuples
$\af=(\af_1,\ldots,\af_k)$ such that $\af_i\geq 1$ for all $i$,
and $\sum\af_i=n$.

Now, $S_k$ acts on $\Z_{> 0}^k$ by permuting the coordinates, and
the size of the orbit of a partition $\ld\in\Z_{> 0}^k$ is
\[
\left.{k}\choose{m_1(\ld),m_2(\ld),\ldots}\right..
\]
Hence,
\[
\sum_{\af}h_{\af}(x)=\sum_{\ld\vdash
n}\left.{k}\choose{m_1(\ld),m_2(\ld),\ldots}\right.h_{\ld}(x).
\]
If we collect the coefficients of $y_j^n$ in the expansion above,
we obtain
\[
\prod_{i}(1-x_iy_j)^{-s}=\sum_{n\geq0}\left(\sum_{\ld\vdash n}
\left.{s}\choose{l(\ld)}\right.\left.{l(\ld)}\choose{m_1(\ld),m_2(\ld),\ldots}\right.
    h_{\ld}(x)\right)y_j^n.
\]
Now, let $\mu=(\mu_1,\mu_2,\ldots)$. Then, by the calculation
above, the coefficient of $y_{j_1}^{\mu_1}y_{j_2}^{\mu_2}\cdots$
in the generating series $\Pi(x,y)^{s}$ is
\[
\prod_{i=1}^{l(\mu)}\left(\sum_{\ld\vdash \mu_i}
\left.{s}\choose{l(\ld)}\right.
\left.{l(\ld)}\choose{m_1(\ld),m_2(\ld),\ldots}\right.
    h_{\ld}(x)\right).
\]
Hence, the Lemma.
\end{pff}

Let $h_n^{(r)}(x)$ be the coefficient of $t^n$ in the generating
series
\[
H(t)^{p^r}=\prod_{i\geq 1}(1-x_it)^{-p^r}
\]
and $h_{\ld}^{(r)}=h_{\ld_1}^{(r)}h_{\ld_2}^{(r)}\cdots$. Then,
\[
\Pi(x,y)^{p^r}=\sum_{\ld}h_{\ld}^{(r)}(x)m_{\ld}(y).
\]
Notice that
\[
\prod_{i\geq 1}(1-x_it)^{-p^r}=\left(\prod_{i\geq
1}(1-x_it)^{-p^{r-1}}\right)^p.
\]
Thus, applying an argument similar to the proof of Lemma
\ref{lem:s-form generating series}, we obtain that
\begin{eqnarray}\label{eqn:higher hom-fcns}
h_n^{(r)}=\sum_{\ld\vdash n}{p\choose l(\ld)}{l(\ld)\choose
m_1(\ld),\ldots,m_n(\ld)}h_{\ld}^{(r-1)}.
\end{eqnarray}
Let
$\Ld^{(r)}=\bigoplus_\ld\Z h_\ld^{(r)}$. Then, we have a sequence
of sublattices
\[
\Ld=\Ld^{(0)}\supset\Ld^{(1)}\supset\cdots\supset\Ld^{(r)}\supset\cdots.
\]

\subsection{Divisibility Properties of Higher Homogeneous Functions}

First, we will examine some divisibility properties of the
coefficient of $h_{\ld}^{(r-1)}(x)m_\mu(y)$ in the generating
series
\[
\Pi(x,y)^{p^r}.
\]
To this end, we have the following Proposition:
\begin{prp}\label{prp:binomial divisibility}
If $\ld$ is a partition of $n$, then
\[
\left.{p}\choose{l(\ld)}\right.
\left.{l(\ld)}\choose{m_1(\ld),m_2(\ld),\ldots}\right.
\]
is divisible by $p$ unless
$\ld=\left(\left(\frac{n}{p}\right)^p\right)$.
\end{prp}

\begin{pff}
First, if $l(\ld)<p$, then $p$ divides
$\left.{p}\choose{l(\ld)}\right.$. On the otherhand, if $l(\ld)=p$
and $\ld\neq\left(\left(\frac{n}{p}\right)^p\right)$, then
$m_i(\ld)<p$ for all $i$, so $p$ divides
\[
\left.{l(\ld)}\choose{m_1(\ld),m_2(\ld),\ldots}\right.=\left.{p}\choose{m_1(\ld),m_2(\ld),\ldots}\right.
.
\]
\end{pff}

This Proposition entitles us to define the following integers:

\begin{dfn}\label{dfn:coefficients} If $\ld\vdash n$ and $\ld\neq((n/p)^p)$, set
\begin{eqnarray}\label{eqn:coefficients}
c_n(\ld)=\frac{1}{p}{p\choose l(\ld)}{l(\ld)\choose
m_1(\ld),\ldots,m_n(\ld)}
\end{eqnarray}
and $c_{pn}((n^p))=0$.
\end{dfn}

\begin{lem}\label{lem:cong1} $c_{pn}(p\ld)=c_n(\ld)$.
\end{lem}

\begin{pff} We have $l(\ld)=l(p\ld)$.
Also $m_{pk}(p\ld)=m_k(\ld)$, and if $(k,p)=1$, then
$m_k(p\ld)=0$. Thus,
\begin{eqnarray*}
c_{pn}(p\ld)&=&\frac{1}{p}{p\choose l(p\ld)}{l(p\ld)\choose
m_p(p\ld),m_{2p}(p\ld)\ldots,m_{pn}(p\ld)}\\
    &=&\frac{1}{p}{p\choose l(\ld)}{l(\ld)\choose
m_1(\ld),\ldots,m_n(\ld)}\\
    &=&c_n(\ld).
\end{eqnarray*}
\end{pff}

\begin{cor}\label{cor:cong2} Let $\ld$ be a partition. Then,
\begin{enumerate}
\item $h_{p\ld}^{(r)}\equiv (h_\ld^{(r-1)})^p$ mod $p\Ld^{(r-1)}$;
\item If $\ld$ contains a part prime to $p$, then
$h_{\ld}^{(r)}\equiv 0$ mod $p\Ld^{(r-1)}$.
\end{enumerate}
\end{cor}

\subsection{A More Suitable Pair of Bases} Our goal is to find the
invariant factors of the matrix whose coefficients are the
coefficients of $h_\ld(x)m_\mu(y)$ in the generating series
\[
\Pi(x,y)^{p^r}=\sum_\ld h_\ld^{(r)}(x)m_\ld(y).
\]
Unfortunately, the transition matrix $M(h^{(r)},h)$ is
complicated. In this subsection, we construct a new family of
bases $(g_{\ld}^{(i,r)})$, and $(M_\ld)$ for $\Ld_\C$ so that
$(g_\ld^{(0,r)})$ and $(M_\ld)$ are bases for $\Ld$,
\[
\Pi(x,y)^{p^r}=\sum_\ld g_\ld^{(r,r)}(x)M_\ld(y)
\]
and the transition matrix $M(g^{(r,r)},g^{(0,r)})$ is
\emph{relatively} simple. In view of Proposition \ref{prp:s-form},
the coefficient of $g_\ld^{(0,r)}(x)M_\mu(y)$ is $\langle
(g_\ld^{(0,r)})^*,M_\mu^*\rangle_{p^r}$. The matrix with these
entries has the same invariant factors as our original matrix.

\begin{lem}\label{lem:main}
For each $r\geq 1$ and $0\leq i\leq r$, there exist multiplicative
bases $(g_{\ld}^{(i,r)})$ for $\Ld^{(i)}$ satisfying
\begin{description}
\item[i] For $(n,p)=1$, $g_n^{(r,r)}=p^{r-i}g_n^{(i,r)}$;
\item[ii] For any $l\geq 1$, $i\geq1$,
$g_{pl}^{(i,r)}=pg_{pl}^{(i-1,r)}+(g_l^{(i-1,r)})^p$;
\item[iii] $g_l^{(i,r)}=\sum_{\ld\vdash l}c_l^{(r-i)}(\ld)h_\ld$ with
$c_l^{(i)}((l))=1$ and $c_{pl}^{(r-i-1)}(p\ld)=c_{l}^{(r-i)}(\ld)$
for all admissible $i$'s.
\end{description}

In particular, $(g_\ld^{(i,r)})$ is a multiplicative basis for
$\Ld^{(i)}$ and the transition matrix $M(g^{(i,r)},h^{(i)})$ is
upper unitriangular.
\end{lem}

\begin{pff} \textbf{The case $r=1$.} Apply induction on the $p$-adic valuation of $l$.
If $l=n$ is an integer prime to $p$, set $g^{(1,1)}_n=h^{(1)}_n$
and
\[
g^{(0,1)}_n=\frac1pg_n^{(1,1)}=\sum_{\ld\vdash n}c_n(\ld)h_\ld
\]
(see equations (\ref{eqn:higher hom-fcns}) and
(\ref{eqn:coefficients})).

For the inductive step, assume that for some $l\geq 1$ we have
found $g^{(0,1)}_l$ satisfying
\[
g_l=\sum_{\ld\vdash l}c_l^{(1)}(\ld)h_{\ld}
\]
for integer coefficients $c_l^{(1)}(\ld)$. Set
\begin{eqnarray}\label{eqn:g_pl}
g_{pl}^{(1,1)}=\sum_{\ld\vdash l}c_l^{(1)}(\ld)h^{(1)}_{p\ld}.
\end{eqnarray}
Then, by Corollary \ref{cor:cong2},
\[
g_{pl}^{(1,1)}\equiv\sum_{\ld\vdash l
}c_l^{(1)}(\ld)(h_{\ld})^p\equiv(g_l^{(0,1)})^p
\]
modulo $p\Ld$. Set
\[
g_{pl}^{(0,1)}=\frac1p[g_{pl}^{(1,1)}-(g_l^{(0,1)})^p],
\]
and note that equation (\ref{eqn:g_pl}) guarantees that
\textbf{iii} holds.

\textbf{The case $r\geq 2$.} Assume by induction that the lemma
holds for all smaller $r$. In particular,  we have integers
$c_l^{(s)}$ for all $s<r$. Hence, for $i>0$, set
\[
g_l^{(i,r)}=\sum_{\ld\vdash l}c_l^{(r-i)}(\ld)h_\ld^{(i)}.
\]
Then, \textbf{iii} holds for $i>0$. Moreover, \textbf{ii} holds
for $i>1$. Indeed, \textbf{ii} is equivalent to
\begin{eqnarray}\label{eqn:induction step}
\hspace{.25in}\sum_{\ld\vdash pl}c_{pl}^{(r-i)}(\ld)h_\ld^{(i)}
    &=&p\left(\sum_{\ld\vdash
        pl}c_{pl}^{(r-i+1)}(\ld)h_\ld^{(i-1)}\right)
+\left(\sum_{\ld\vdash l}c_l^{(r-i+1)}(\ld)h_\ld^{(i-1)}\right)^p.
\end{eqnarray}
By induction, we have
\[
g_{pl}^{(i-1,r-1)}=pg_{pl}^{(i-2,r-1)}+\left(g_l^{(i-2,r-1)}\right)^p
\]
or
\[
\sum_{\ld\vdash pl}c_{pl}^{(r-i)}(\ld)h_\ld^{(i-1)}
=p\left(\sum_{\ld\vdash
pl}c_{pl}^{(r-i+1)}(\ld)h_\ld^{(i-2)}\right)
+\left(\sum_{\ld\vdash l}c_l^{(r-i+1)}(\ld)h_\ld^{(i-2)}\right)^p.
\]
Since $M(h^{(i)},h^{(i-1)})=M(h^{(i-1)},h^{(i-2)})$, this implies
equation (\ref{eqn:induction step}). Similarly, one checks that
\textbf{i} holds for all $i>0$.

Now, we define $g_l^{(0,r)}$ inductively by the formulas
\begin{eqnarray*}
g_{n}^{(0,r)}&=&\frac1pg_n^{(1,r)}\;\;\;\mbox{if }(n,p)=1,\\
g_{pl}^{(0,r)}&=&\frac1p\left[g_{pl}^{(1,r)}-g_l^{(0,r)}\right].
\end{eqnarray*}

To finish the proof, we check by induction on the $p$-adic
valuation of $l$ that $g_l^{(0,r)}\in\Ld$, the coefficient of
$h_l$ in $g_l^{(0,r)}$ is 1, and
$c_{pl}^{(r-1)}(p\ld)=c_l^{(r)}(\ld)$. For the induction base,
observe that when $(n,p)=1$, $g_n^{(1,r)}$ is a linear combination
of $h_\ld^{(1)}$ satisfying the conditions of Corollary
\ref{cor:cong2}(2), so $g_n^{(0,r)}\in\Ld$. It is also clear from
equation (\ref{eqn:higher hom-fcns}) that the coefficient of $h_n$
in $g_n^{(0,r)}$ is 1. Next, we prove that
\begin{eqnarray}\label{eqn:pn-n coeff}
c_{pn}^{(r-1)}(p\ld)=c_n^{(r)}(\ld).
\end{eqnarray}
Indeed, the number $c_n^{(r)}(\ld)$ is the coefficient of $h_\ld$
in
\begin{eqnarray}\label{eqn:coef, base case}
g_n^{(0,r)}=\frac1pg_n^{(1,r)}=\frac1p\left[\sum_{\mu\vdash
n}c_n^{(r-1)}(\mu)h_\mu^{(1)} \right].
\end{eqnarray}
On the other hand, $c_{pn}^{(r-1)}(p\ld)$ is the coefficient of
$h_{p\ld}$ in
\begin{eqnarray}\label{eqn:coef, base case 2}
g_{pn}^{(0,r)}&=&\frac1p\left[g_{pn}^{(1,r)}-g_n^{(0,r)}\right]\\
\nonumber    &=&\frac1p\left[\sum_{\nu\vdash pn}
c_{pn}^{(r-1)}(\nu)h_\nu
    -\left(\sum_{\mu\vdash n} c_n^{(r-1)}(\mu)h_\mu^{(1)}
    \right)^p\;\right].
\end{eqnarray}
No term of the form $h_{p\ld}$ comes from
\[
\left(\sum_{\mu\vdash n} c_n^{(r-1)}(\mu)h_\mu^{(1)} \right)^p,
\]
as every partition occurring there contains a part prime to $p$.
Hence, $c_{pn}^{(r-1)}(p\ld)$ is the coefficient of $h_{p\ld}$ in
\[
\frac1p\sum_{\nu\vdash pn} c_{pn}^{(r-1)}(\nu)h_\nu.
\]
It follows from equation (\ref{eqn:higher hom-fcns}) that
$h_{p\ld}$ appears in $h_\nu^{(1)}$ only if $\nu$ is of the form
$p\mu$ for $\mu\vdash n$. So (\ref{eqn:pn-n coeff}) follows by
comparing (\ref{eqn:coef, base case}) and (\ref{eqn:coef, base
case 2}), since by induction we know that
$c_{pn}^{(r-2)}(p\mu)=c_n^{(r-1)}(\mu)$.

For the inductive step, we have
\begin{eqnarray*}
g_{pl}^{(1,r)}&=&\sum_{\ld\vdash pl}c_{pl}^{(r-1)}(\ld)h_\ld^{(1)}\\
        &\equiv&\sum_{\ld\vdash l}c_l^{(r)}(\ld)(h_\ld)^p\\
        &\equiv& (g_l^{(0,r)})^p
\end{eqnarray*}
modulo $p\Ld$. This shows that $g_{pl}^{(0,r)}\in\Ld$. Moreover,
equation (\ref{eqn:higher hom-fcns}) together with the definition
of $g_{pl}^{(0,r)}$ imply that the coefficient of $h_{pl}$ in
$g_{pl}^{(0,r)}$ is 1.

Finally, we show that $c_{p^2l}^{(r-1)}(p\ld)=c_{pl}^{(r)}(\ld)$.
The number $c_{p^2l}^{(r-1)}(p\ld)$ is the coefficient of
$h_{p\ld}$ in
\[
g_{p^2l}^{(0,r-1)}=\frac1p[g_{p^2l}^{(1,r-1)}-(g_{pl}^{(0,r-1)})^p]\\
\]
and the number $c_{pl}^{(r)}(\ld)$ is the coefficient of $h_\ld$
in
\[
g_{pl}^{(0,r)}=\frac1p[g_{pl}^{(1,r)}-(g_l^{(0,r)})^p]
\]
By induction, we may assume that the coefficient of $h_{p\ld}$ in
$(g_{pl}^{(0,r-1)})^p$ is the same as the coefficient of $h_{\ld}$
in $(g_l^{(0,r)})^p$. Now, argue as above that the coefficient of
$h_{p\ld}$ in $g_{p^2l}^{(1,r-1)}$ is the same as the coefficient
of $h_\ld$ in $g_{pl}^{(1,r)}$.
\end{pff}

\begin{cor}\label{cor:good gen series} There exists a basis $(M_{\ld})$ for $\Ld$ such that
\[
\Pi(x,y)^{p^r}=\sum_{\ld}g_{\ld}^{(r,r)}(x)M_\ld(y)
\]
\end{cor}

\begin{pff}
Let $Z=M(h^{(r)},g^{(r,r)})$. By Lemma \ref{lem:main}
\textbf{iii}, it follows that $Z$ is the same for all $r$. Now,
define $(M_{\ld})$ by $M(M,m)=Z^t$. Then, \beqr
\Pi(x,y)^{p^r}&=&\sum_{\ld}h^{(r)}_{\ld}(x)m_{\ld}(y)\\
    &=&\sum_{\ld,\mu}Z_{\ld\mu}g^{(r,r)}_{\mu}(x)m_{\ld}(y)\\
    &=&\sum_{\mu}g^{(r,r)}_{\mu}(x)\left(\sum_{\ld}Z_{\ld\mu}m_{\ld}(y)\right)\\
    &=&\sum_{\mu}g^{(r,r)}_{\mu}(x)M_{\mu}(y)
\eeqr
\end{pff}

\subsection{Invariant factors of the $p$-form on $\Ld$}\label{subsection:p-form} For ease
of notation, we omit superscripts in this section. That is, set
$g_\ld:=g_\ld^{(0,1)}$. Then, by Corollary \ref{cor:good gen
series} the entries of the Gram matrix of the $p$-form on $\Ld$
are the coefficients of $g_\ld(x)M_\mu(y)$ in
\[
\Pi(x,y)^{p}=\sum_{\mu}\prod_{(n,p)=1}
\left(pg_n\right)^{m_n(\mu)}\prod_{i\geq1}\left(pg_{p^in}+(g_{p^{i-1}n})^p\right)^{m_{p^in}(\mu)}M_{\mu}(y).
\]

Let $n$ be an integer prime to $p$, and $i\geq1$. Define elements
$G_{(n)}=g_n$, $G_{(n^{p^i})}=g_{p^in}$, and
\begin{eqnarray}\label{eqn:relation}
    G_{(p^in)}=pG_{(n^{p^i})}+(G_{(n^{p^{i-1}})})^p=pg_{p^in}+(g_{(p^{i-1}n)})^p.
\end{eqnarray}
If $m$ is any integer, let $m=\sum_ia_ip^i$ be its $p$-adic
decomposition. We set
\begin{eqnarray}\label{eqn:relation 2}
G_{(n^m)}=\prod_i(G_{(n^{p^i})})^{a_i}\andeqn
G_{((p^in)^m)}=(G_{(p^in)})^m.
\end{eqnarray}
For $\ld=(1^{m_1}2^{m_2}\cdots)$, set
\[
G_{\ld}=G_{(1^{m_1})}G_{(2^{m_2})}\cdots.
\]
We will demonstrate below that $(G_\ld)$ forms a basis for $\Ld$.

\begin{lem}\label{lem:lead coefficient} There exist integers
$C_i(\ld)$ such that
\[
(G_{(n^{p^{i-j}})})^{p^j}=\sum_{\ld\vdash p^in}C_{ij}(\ld)G_\ld
\]
and
\begin{enumerate}
\item $C_{ij}((n^{p^i}))=\pm p^{\frac{p^j-1}{p-1}}$, and
\item $C_{ij}(\ld)=0$ unless $\ld=(p^{i_1}n\geq p^{i_2}n\geq\cdots)$.
\end{enumerate}
\end{lem}

\begin{pff}
Proceed by induction on $i\geq 0$ and $0\leq j\leq i$. For $i=0$,
or $i>0$ and $j=0$ the lemma holds trivially. Suppose that $i,j>0$
and write
\[
\ph_{i-1,j-1}(G)=\sum_{\ld\neq(n^{p^{i-1}})}C_{i-1,j-1}(\ld)G_\ld.
\]
Now,
\begin{eqnarray}\label{eqn:binomial expansion}
\nonumber(G_{(n^{p^{i-j}})})^{p^j}&=&((G_{(n^{p^{i-j}})})^{p^{j-1}})^p\\
\nonumber           &=&\left(\ph_{i-1,j-1}(G)\pm p^{(p^{j-1}-1)/(p-1)}G_{(n^{p^{i-1}})}\right)^p\\
            &=&\sum_{k=0}^{p-1}\left.{p}\choose{k}\right.
                (\pm p^{(p^{j-1}-1)/(p-1)})^{k}\ph_{i-1,j-1}(G)^{p-k}G_{(n^{p^{i-1}})}^k\\
\nonumber                &&\hspace{.25in}\pm
                p^{(p^j-p)/(p-1)}(G_{(n^{p^{i-1}})})^p.
\end{eqnarray}
Induction applies to the first summand in (\ref{eqn:binomial
expansion}). By (\ref{eqn:relation}), we have that
\[
(G_{(n^{p^{i-1}})})^p=G_{(p^in)}-pG_{(n^{p^i})}.
\]
Therefore,
\begin{eqnarray*}
p^{(p^{j}-p)/(p-1)}(G_{(n^{p^{i-1}})})^p&=&p^{(p^{j}-p)/(p-1)}(G_{(p^in)}-pG_{(n^{p^i})})\\
                                &=&p^{(p^{j}-p)/(p-1)}G_{(p^in)}-p^{(p^j-p)/(p-1)+1}G_{(n^{p^i})}.
\end{eqnarray*}
Hence, the lemma.
\end{pff}

\begin{cor} $(G_\ld)$ forms a basis for $\Ld$.
\end{cor}

\begin{pff}
It is easy to write $G_\ld$ as an integral linear combination of
$g_\mu$'s. Conversely, write
\begin{eqnarray}\label{eqn:M(g,G)}
\nonumber g_\ld&=&\prod_{(n,p)=1}\prod_{i\geq
0}(g_{p^in})^{m_{p^in}(\ld)}\\
    &=&\prod_{(n,p)=1}\prod_{i\geq
0}(G_{(n^{p^i})})^{m_{p^in}(\ld)}.
\end{eqnarray}
and expand using (\ref{eqn:relation}) and (\ref{eqn:relation 2})
above. By Lemma \ref{lem:lead coefficient}(2), it is enough to
expand a product of the form
\[
\prod_{i\geq 0}(G_{(n^{p^i})})^{m_i},
\]
which can be obtained by induction on $\sum_im_i$. Indeed, if all
$m_i<p$, there is nothing to do. Otherwise, find the smallest $i$
such that $m_i>p$, and let $m_i=\sum_{j\geq 0}a_j^ip^j$ be the
$p$-adic expansion of $m_i$. We have
\begin{eqnarray*}
\left(G_{(n^{p^i})}\right)^{m_i}&=&\prod_{j\geq0}\left(G_{(n^{p^i})}\right)^{p^ja_j^i}\\
        &=&\prod_{j\geq
        0}\left((\mathrm{const.})G_{(n^{p^{i+j}})}\right)^{a_j^i}+(*).
\end{eqnarray*}
It follows from Lemma \ref{lem:lead coefficient}(2) that $(*)$ is
a linear combination of terms having fewer than $m_i$ components
of the form $G_{(n^{p^k})}$. Observing that $\sum_{j\geq
0}a_j^i<m_i$, completes the induction.
\end{pff}


In this basis, we have that
\begin{eqnarray}\label{eqn:presentation}
\hspace{.25in}\Pi(x,y)^{p}
    &=&\sum_{\mu}\prod_{(n,p)=1}
p^{m_n(\mu)}(G_{(n)}(x))^{m_n(\mu)}\prod_{i\geq1}(G_{((p^in)^{m_{p^in}(\mu)})}(x))M_{\mu}(y).
\end{eqnarray}

We need to analyze the term $(G_{(n)})^{m_n(\mu)}$ in the
generating series above. To this end, let
$m_n(\mu)=m_n=\sum_ia^n_ip^i$ be its p-adic decomposition. Then,
\[
(G_{(n)})^{m_n}=\prod_i(G_{(n)})^{a^n_ip^i}.
\]

Recall the definition of $d_p(a)$ given by equation (\ref{formula:
geometric sum}) on page \pageref{formula: geometric sum}.

\begin{cor}\label{cor:lead coefficient} The coefficient of $G_{(n^{m_n})}$ in the product
$(G_{(n)})^{m_n}$ is $p^{d_p(m_n(\ld))}$.
\end{cor}

\begin{pff} By Lemma \ref{lem:lead coefficient} the coefficient of $G_{(n^{m_n})}$ in the product
$(G_{(n)})^{m_n}$ is
\[
\pm\prod_{j\geq1}p^{a_j^n\left(\frac{p^j-1}{p-1}\right)}.
\]
Observe that
\begin{eqnarray*}
\sum_{j\geq1}a_j^n\left(\frac{p^j-1}{p-1}\right)&=&\sum_{j\geq1}a_j^n(1+p+\cdots+p^{j-1})\\
                &=&\sum_{j\geq1}\sum_{k\geq1}\left\lfloor\frac{a_j^np^j}{p^k}\right\rfloor\\
                &=&\sum_{k\geq1}\left\lfloor\frac{\sum_{j\geq1}a_j^np^j}{p^k}\right\rfloor\\
                &=&d_p(m_n).
\end{eqnarray*}
\end{pff}

\begin{cor} The matrix whose entries are the coefficients of
$G_{\ld}(x)M_{\mu}(y)$ in the product $\Pi(x,y)^{p}$ is upper
triangular. Moreover, if $\ld=(1^{m_1}2^{m_2}\cdots)$ and
$m_n=\sum_ia^n_ip^i$ is the p-adic decomposition of $m_n$, then
the coefficient of $G_{\ld}(x)M_{\ld}(y)$ is
\[
\pm\prod_{(n,p)=1}p^{m_n+d_p(m_n)}.
\]
\end{cor}

The following proposition will complete the proof of Theorem
\ref{thm:main} for the special case $r=1$.

\begin{prp}\label{prp:divisibility along rows} The coefficient of $G_{\ld}(x)M_{\mu}(y)$ in the
product $\Pi(x,y)^p$ is divisible by that of
$G_{\ld}(x)M_{\ld}(y)$.
\end{prp}

\begin{pff} First, observe that by Corollary \ref{cor:lead coefficient},
the coefficient of $G_{\ld}(x)M_{\ld}(y)$ is
\[
\prod_{(n,p)=1}p^{m_n(\ld)+d_p(m_n(\ld))}.
\]

Let $m_n(\mu)=\sum_{j\geq0}a_j^n(\mu)p^j$ be the $p$-adic
expansion of $m_n(\mu)$. Then, (\ref{eqn:presentation}) and Lemma
\ref{lem:lead coefficient} give us that the coefficient of
$M_\mu(y)$ in $\Pi(x,y)^p$ is
\[
\prod_{(n,p)=1}p^{m_n(\mu)}\prod_{j\geq 0}\left(\sum_{\sm\vdash
p^jn}C_j(\sm)G_\sm(x)\right)^{a_j^n(\mu)}\prod_{i\geq
1}\left(G_{\left((p^in)^{m_{p^in}(\mu)}\right)}(x)\right).
\]
Using property (2) of the lemma, we see that if $C_j(\sm)\neq0$
and $\sm\neq (n^{p^j})$, then $\sm$ must have fewer than $p^j$
parts prime to $p$. Hence, the $\ld$ that occur when expanding the
expression above have fewer parts prime to $p$ than $\mu$ and,
therefore, the coefficient of $G_{\ld}(x)M_{\mu}(y)$ is divisible
by $\prod_{(n,p)=1}p^{m_n(\ld)}$. Finally, property (1) of the
lemma implies that each time $(n^{p^j})$ appears in $\ld$, the
coefficient of $G_{\ld}(x)M_{\mu}(y)$ is divisible by
$p^{d_p(p^i)}$.
\end{pff}

The following corollary is immediate.

\begin{cor}
There exists a basis $(N_{\ld})$ for $\Ld$ such that
\[
\Pi(x,y)^{p}=\sum_{\ld}\prod_{(n,p)=1}p^{m_n(\ld)+d_p(m_n(\ld))}G_{\ld}(x)N_{\ld}(y).
\]
\end{cor}

\subsection{Outline of the Approach to the $p^r$-form}\label{subsection:strategy} The basis
$(G_\ld)$ for $\Ld$ constructed in the previous section (together
with computational evidence) suggests an approach to finding the
invariant factors of the $p^r$-form on $\Ld$. Indeed, suppose that
we can construct a basis $(G_\ld^{(r)})$ with the following
properties:

\begin{enumerate}
\item[(P1)] For $i\geq0$,
\[
G_{(p^in)}^{(r)}=\frac1{\lceil{p^{r-i}}\rceil} g_{p^in}^{(r,r)}
\]
where $\lceil x\rceil$ is the ceiling function.
\item[(P2)] For $0\leq i<r$ and $j\geq 0$,
\[
G_{((p^{i}n)^{p^j})}^{(r)}=\sum_{\ld\vdash
p^{i+j}n}X_{ij}(\ld)g^{(0,r)}_\ld
\]
for integer coefficients $X_{ij}(\ld)$ satisfying
$X^n_{ij}((p^{i+j}n))=p^i$ and $X_{ij}(\ld)=0$ unless
$\ld=(p^{k_1}n\geq p^{k_2}n\geq\cdots)$.
\item[(P3)] For $1\leq i<r$ and $j\geq 0$,
\[
G_{((p^{i}n)^{p^j})}^{(r)}=(G^{(r)}_{(p^{i-1}n)^{p^j}})^p+p\,G^{(r)}_{(p^{i-1}n)^{p^{j+1}}}.
\]
\item[(P4)] For $i\geq r$,
\[
G^{(r)}_{(p^in)}=(G^{(r)}_{((p^{r-1}n)^{p^{i-r}})})^p+p\,G^{(r)}_{((p^{r-1}n)^{p^{i-r+1}})}.
\]
\item[(P5)]\label{Basis condition} For a partition
$\ld=(1^{m_1}2^{m_2}\ldots)$,
\[
G_\ld^{(r)}=\prod_i G_{(i^{m_i})}^{(r)}
\]
where, if $m=\sum_{j\geq 0}a_jp^j$ is the $p$-adic expansion of
$m$, then
\[
G_{(p^in)^m}^{(r)}=\left\{\begin{array}{ll}\prod_{j\geq0}(G^{(r)}_{(p^in)^{p^j}})^{a_j}&
                    \mbox{if }0\leq i<r;\\
                   (G_{(p^in)}^{(r)})^m&\mbox{otherwise.}\end{array}\right.
\]
\end{enumerate}
By taking $i=0$ in property (P2) it follows that we indeed have a
basis for $\Ld$. Indeed, it is easy to write the $G_\ld^{(r)}$ in
terms of $g_\ld^{(0,r)}$. To write the $g_\ld^{(0,r)}$ in terms of
$G_\ld^{(r)}$, we show how to write $g_{p^in}^{(0,r)}$ in terms of
$G_\ld^{(r)}$ by induction on $i$, starting from
$g_n^{(0,r)}=G_{(n)}^{(r)}$. Now, assume that we can write
$g_{p^jn}^{(0,r)}$ in terms of $G_\ld^{(r)}$ for all $j<i$. Then,
\[
G_{(n^{p^i})}^{(r)}=\sum_{\ld\vdash
p^{i}n}X_{0i}(\ld)g^{(0,r)}_\ld
\]
with $X_{0i}((p^in))=1$ and $X_{0i}(\ld)=0$ unless
$\ld=(p^{k_1}n\geq p^{k_2}n\geq\cdots)$. In particular, induction
applies to the $g_\ld^{(0,r)}$ occurring with nonzero
coefficients. Since both sets $(G_{\ld}^{(r)})$ and
$(g_\ld^{(0,r)})$ are labelled by partitions, we conclude that we
have a basis.

In this basis, we have
\[
\Pi(x,y)^{p^r}=\sum_{\ld}\prod_{(n,p)=1}\prod_{i=0}^{r-1}\left(p^{r-i}G_{(p^in)}^{(r)}(x)\right)^{m_{p^in}(\ld)}
    \prod_{l\geq1}\left(G_{\left((p^rl)^{m_{p^rl}(\ld)}\right)}^{(r)}(x)\right)M_{\ld}(y).
\]

It remains to multiply out the terms
$(p^{r-i}G_{(p^in)}^{(r)})^{m_{p^in}(\ld)}$ using (P3) and (P4).
Indeed, we have

\begin{lem}\label{lem:p^r lead coefficient} For $0\leq i<r$, there exist integers $C_{ijk}(\ld)$ such that
\[
(G^{(r)}_{((p^in)^{p^{j-k}})})^{p^k}=\sum_{\ld\vdash
p^{i+j}n}C_{ijk}(\ld)G^{(r)}_\ld
\]
and
\begin{enumerate}
\item $C_{ijk}(((p^in)^{p^j}))=\pm p^{\frac{p^k-1}{p-1}}$, and
\item $C_{ijk}(\ld)=0$ unless $\ld=(p^{i_1}n\geq
p^{i_2}n\geq\cdots)\geq((p^in)^{p^j})$.
\end{enumerate}
\end{lem}

The proof of this Lemma is similar to the proof of Lemma
\ref{lem:lead coefficient} using (P3) and (P4) above. The
analogous corollaries to Lemma \ref{lem:lead coefficient} are also
true, and their proofs require no new techniques. In particular,
one has

\begin{cor}\label{cor:p^r-lead coefficient} For $0\leq i<r$, the coefficient of $G_{(p^in)^{m}}^{(r)}$ in the product
$(G_{(p^in)}^{(r)})^{m}$ is $p^{d_p(m)}$.
\end{cor}

\begin{cor} The matrix whose entries are the coefficients of
$G_{\ld}^{(r)}(x)M_{\mu}(y)$ in the product $\Pi(x,y)^{p^r}$ is
upper triangular. Moreover, if $\ld=(1^{m_1}2^{m_2}\cdots)$, then
the coefficient of $G_{\ld}^{(r)}(x)M_{\ld}(y)$ is
\[
D_r(\ld)=\prod_{(n,p)=1}\prod_{i=0}^{r-1}p^{[(r-i)m_{p^in}+d_p(m_{p^in})]}.
\]
(recall $D_r(\ld)$ from (\ref{eqn: invariant factors}) on page
\pageref{eqn: invariant factors}.)
\end{cor}

\begin{prp}\label{prp:p^r-divisibility along rows} The coefficient of $G_{\ld}^{(r)}(x)M_{\mu}(y)$ in the
product $\Pi(x,y)^{p^r}$ is divisible by the coefficient of
$G_{\ld}^{(r)}(x)M_{\ld}(y)$.
\end{prp}

Hence, we arrive at the conjecture

\begin{cnj}\label{cnj:p^r-invariant factors} The invariant factors of the Gram matrix of the $p^r$-form on
the degree $d$ component of $\Ld$ are
\[
\{D_r(\ld)|\ld\vdash d\}.
\]
\end{cnj}

In the following sections, we show that the conjecture is true
under the assumption that $r\leq p$.

\begin{rmk} For the sake of the lemmas and propositions above, conditions
(P3) and (P4) are overly strict. Instead, we could replace these
with the following:
\begin{enumerate}
\item[(P3')] For $1\leq i<r$ and $j\geq 0$,
\[
\left(G_{((p^{i-1}n)^{p^j})}^{(r)}\right)^p=\sum_{\ld\vdash
p^{i+j}n}Y_{ij}(\ld)G_\ld^{(r)}
\]
where
\begin{itemize}
\item $Y_{ij}(((p^{i-1}n)^{p^{j+1}}))=\pm p$;
\item $Y_{ij}(\ld)=0$ unless $\ld=(p^{i_1}n\geq
p^{i_2}n\geq\cdots)\geq((p^{i-1}n)^{p^{j+1}})$.
\end{itemize}
\item[(P4')] For $i\geq r$,
\[
\left(G_{((p^{r-1}n)^{p^{i-r}})}^{(r)}\right)^p=\sum_{\ld\vdash
p^in}Y_{i,i-r}(\ld)G_\ld^{(0,r)}
\]
where
\begin{itemize}
\item $Y_{i,i-r}(((p^{r-1}n)^{p^{i-r+1}}))=\pm p$;
\item $Y_{i,i-r}(\ld)=0$ unless $\ld=(p^{i_1}n\geq
p^{i_2}n\geq\cdots)\geq((p^{r-1}n)^{p^{i-r+1}})$.
\end{itemize}
\end{enumerate}
In fact, for the case $r=p$, we will only obtain conditions (P3')
and (P4'). The problem with these conditions is that they do not
make satisfactory inductive hypotheses for calculations below.
\end{rmk}

\subsection{Invariant factors of the $p^r$-form on $\Ld$, for
$r<p$} Fix $r\geq 2$, and assume that for $s<r$ we have
constructed a basis $(G_\ld^{(s)})$ for $\Ld$ satisfying
properties (P1)-(P5) of \S \ref{subsection:strategy}. We take as
our base case the basis $(G_\ld^{(1)}):=(G_\ld)$ constructed in \S
\ref{subsection:p-form}.

Define a basis $(\hG_\ld)$ for $\Ld$ inductively by the formula
\[
M(g^{(0,r)},\hG)=M(g^{(0,r-1)},G^{(r-1)}).
\]
Then, since
$M(g^{(r-1,r)},g^{(0,r)})=M(g^{(r-1,r-1)},g^{(0,r-1)})$, property
(P1) implies that for $1\leq i<r$,
\begin{eqnarray*}\label{eqn:inductive step 1}
g_{p^in}^{(r,r)}&=&pg_{p^in}^{(r-1,r)}+(g_{p^{i-1}n}^{(r-1,r)})^p\\
\nonumber    &=&p\left(p^{r-i-1}\hG_{(p^in)}\right)+\left(p^{r-i}\hG_{(p^{i-1}n)}\right)^p\\
\nonumber   &\equiv&0
\end{eqnarray*}
modulo $p^{r-i}\Ld$. We may therefore define
\begin{equation*}\label{eqn:inductive step 2}
G_{(p^in)}^{(r)}=\frac1{\lceil p^{r-i}\rceil}g_{p^in}^{(r,r)}
\end{equation*}
for all $i\geq 0$. (Note that $\lceil p^{r-i}\rceil=1$ for $i\geq
r$.)

Next, we construct $G_{((p^{i-j}n)^{p^j})}^{(r)}$ by induction on
$0\leq i<r$ and $0\leq j\leq i$ so that we obtain properties (P2)
and (P3). Notice here that $G_{(n)}^{(r)}=g_n^{(0,r)}$, so (P2)
holds, and (P3) is vacuous.

Now, let $1<i\leq r-1$. Assume that for $k<i$ and $1\leq l\leq k$
we have constructed $G_{((p^{k-l}n)^{p^l})}^{(r)}$ such that
\begin{equation}\label{eqn:inductive step 3}
G_{((p^{k-l}n)^{p^l})}^{(r)}=\hG_{((p^{k-l}n)^{p^l})}+p^{k-l+1}\ph_{k,l}(g)
\end{equation}
where, here and throughout the paper, $\ph_{k,l}(g)$ is a
polynomial in the $g_{p^mn}^{(0,r)}$ with $m<k$. (Notice that
equation (\ref{eqn:inductive step 3}) and induction on $r$ imply
that properties (P2) and (P3) hold for $k<i$.)

By property (P3),
\begin{eqnarray*}
G_{(p^in)}^{(r)}&=&\hG_{(p^in)}+p^{(p-1)(r-i)}\left(\hG_{(p^{i-1}n)}\right)^p\\
        &=&\left[p\hG_{((p^{i-1}n)^p)}+(\hG_{(p^{i-1}n)})^p\right]+p^{(p-1)(r-i)}\left(\hG_{(p^{i-1}n)}\right)^p,
\end{eqnarray*}
and, by equation (\ref{eqn:inductive step 3}),
\begin{equation*}
\hG_{(p^{i-1}n)}=G_{(p^{i-1}n)}^{(r)}-p^{i}\ph_{i-1,0}(g).
\end{equation*}
Thus,
\begin{eqnarray}\label{eqn:verifies polynomial}
\hspace{.25in}G_{(p^in)}^{(r)}&=&p\hG_{((p^{i-1}n)^p)}+\left(G_{(p^{i-1}n)}^{(r)}-p^{i}\ph_{i-1,0}(g)\right)^p\\
    \nonumber&&+p^{(p-1)(r-i)}\left(\hG_{(p^{i-1}n)}\right)^p.
\end{eqnarray}
Since $p> r$ and $i\leq r-1$, it follows that
\[
(p-1)(r-i)>(r-1)(r-(r-1))=r-1\geq i.
\]
Thus,
\begin{equation*}\label{eqn:inductive step 4}
G_{(p^in)}^{(r)}\equiv
p\hG_{((p^{i-1}n)^p)}+(G_{(p^{i-1}n)}^{(r)})^p
\end{equation*}
modulo $p^{i+1}\Ld$. Set
\begin{eqnarray*}
G_{((p^{i-1}n)^p)}^{(r)}&=&\frac1p\left[G_{(p^in)}^{(r)}-(G_{(p^{i-1}n)}^{(r)})^p\right].
\end{eqnarray*}
By induction on $i$, property (P2) holds for
$G_{(p^{i-1}n)}^{(r)}$. Together with equation (\ref{eqn:verifies
polynomial}), this implies that
\[
G_{((p^{i-1}n)^p)}^{(r)}=\hG_{((p^{i-1}n)^p)}+p^i\ph_{i,1}(g).
\]

Next, let $1\leq j<i$ and assume by induction that we have
constructed
\begin{equation*}
G^{(r)}_{((p^{i-j}n)^{p^j})}=\hG_{((p^{i-j}n)^{p^j})}+p^{i-j+1}\ph_{i,j}(g).
\end{equation*}
Then, by property (P3) and equation (\ref{eqn:inductive step 3}),
\begin{eqnarray}\label{eqn:inductive step 5}
G^{(r)}_{((p^{i-j}n)^{p^j})}&=&\left[p\hG_{((p^{i-j-1}n)^{p^{j+1}})}+(G^{(r)}_{((p^{i-j-1}n)^{p^j})}-p^{i-j}\ph_{i-1,j}(g))^p\right]\\\nonumber
                                                &&+p^{i-j+1}\ph_{i,j}(g).
\end{eqnarray}
Thus,
\begin{eqnarray*}
G^{(r)}_{((p^{i-j}n)^{p^j})}&\equiv&p\hG_{((p^{i-j-1}n)^{p^{j+1}})}+(G^{(r)}_{((p^{i-j-1}n)^{p^j})})^p
\end{eqnarray*}
modulo $p^{i-j+1}\Ld$.

Set
\begin{eqnarray*}
G^{(r)}_{((p^{i-j-1}n)^{p^{j+1}})}&=&\frac1p\left[G^{(r)}_{((p^{i-j}n)^{p^j})}-(G^{(r)}_{((p^{i-j-1}n)^{p^j})})^p\right].
\end{eqnarray*}
By induction on $i$, property (P2) holds for
$G^{(r)}_{((p^{i-j-1}n)^{p^j})}$. Therefore, equation
(\ref{eqn:inductive step 5}) implies that
\[
G^{(r)}_{((p^{i-j-1}n)^{p^{j+1}})}=\hG_{((p^{i-j-1}n)^{p^{j+1}})}+p^{i+j}\ph_{i,j+1}(g).
\]
This completes the construction of the
$G^{(r)}_{((p^{i-j}n)^{p^j})}$ for $i<r$.

The last step in the construction is to obtain the elements
$G_{(p^in)}^{(r)}$ for $i\geq r$ satisfying properties (P3) and
(P4) of \S \ref{subsection:strategy}. To this end, assume that
$i\geq r$ and that
\[
G_{(p^{i-1}n)}^{(r)}=\hG_{((p^{r-1}n)^{p^{i-r}})}+p^r\ph_{i-1,0}(g).
\]
(Note that we have the right to make this assumption because when
$i=r$, we have
\[
G_{(p^{r-1}n)}^{(r)}=\hG_{(p^{r-1}n)}+p^{p-1}(\hG_{(p^{r-2}n)})^p
\]
and $p>r$.)

Then,
\begin{eqnarray*}
G_{(p^in)}^{(r)}&=&g_{p^in}^{(r)}\\
                &=&pg_{p^in}^{(r-1,r)}+(g_{p^{i-1}n}^{(r-1,r)})^p\\
                &=&p\hG_{(p^in)}+\left(G_{((p^{r-1}n)^{p^{i-r}})}^{(r)}-p^{r}\ph_{i-1,0}(g)\right)^p\\
                &\equiv&p\hG_{(p^in)}+\left(G_{((p^{r-1}n)^{p^{i-r}})}^{(r)}\right)^p
\end{eqnarray*}
modulo $p^{r+1}\Ld$. Set
\begin{eqnarray*}
G_{((p^{r-1}n)^{p^{i-r+1}})}^{(r)}&=&\frac1p\left[G_{(p^in)}^{(r)}-\left(G_{((p^{r-1}n)^{p^{i-r}})}^{(r)}\right)^p\right]\\
        &=&\hG_{(p^in)}+p^{r}\ph_{i,1}(g).
\end{eqnarray*}
By induction on $r$,
\[
\hG_{(p^in)}=p\hG_{((p^{r-2}n)^{p^{i-r+1}})}+\left(\hG_{((p^{r-2}n)^{p^{i-r}})}\right)^p
\]
and, by induction on $i$,
\[
G_{((p^{r-2}n)^{p^{i-r}})}^{(r)}=\hG_{((p^{r-2}n)^{p^{i-r}})}+p^{r-1}\ph_{i-1,1}(g).
\]
Thus,
\begin{eqnarray*}
G_{((p^{r-1}n)^{p^{i-r+1}})}^{(r)}&=&p\hG_{((p^{r-2}n)^{p^{i-r+1}})}+\left(G_{((p^{r-2}n)^{p^{i-r}})}^{(r)}-p^{r-1}\ph_{i-1,1}(g)\right)^p\\
            &&+p^{r}\ph_{i,1}(g)\\
            &\equiv&\left(G_{((p^{r-2}n)^{p^{i-r}})}^{(r)}\right)^p
\end{eqnarray*}
modulo $p^{r}\Ld$. Set
\begin{eqnarray*}
G_{((p^{r-2}n)^{p^{i-r+2}})}^{(r)}&=&\frac1p\left[G_{((p^{r-1}n)^{p^{i-r+1}})}^{(r)}-\left(G_{((p^{r-2}n)^{p^{i-r}})}^{(r)}\right)^p\right]\\
    &=&\hG_{((p^{r-2}n)^{p^{i-r+1}})}+p^{r-1}\ph_{r,2}(g).
\end{eqnarray*}
Now, assume that $2\leq j<r$ and that
\[
G_{((p^{r-j}n)^{p^{i-r+j}})}^{(r)}=\hG_{((p^{r-j}n)^{p^{i-r+j-1}})}+p^{r-j+1}\ph_{i,j}(g).
\]
By induction on $r$
\[
\hG_{((p^{r-j}n)^{p^{i-r+j-1}})}=p\hG_{((p^{r-j-1}n)^{p^{i-r+j}})}+(\hG_{((p^{r-j-1}n)^{p^{i-r+j-1}})})^p
\]
and, by induction on $i$,
\[
G_{((p^{r-j-1}n)^{p^{i-r+j}})}^{(r)}=\hG_{((p^{r-j-1}n)^{p^{i-r+j-1}})}+p^{r-j}\ph_{i-1,j}(g).
\]
Hence,
\begin{eqnarray*}
G_{((p^{r-j}n)^{p^{i-r+j}})}^{(r)}&=&p\hG_{((p^{r-j-1}n)^{p^{i-r+j}})}
+(G_{((p^{r-j-1}n)^{p^{i-r+j}})}^{(r)}-p^{r-j}\ph_{i-1,j}(g))^p\\
&&+p^{r-j+1}\ph_{i,j}(g)\\
    &\equiv&(G_{((p^{r-j-1}n)^{p^{i-r+j}})}^{(r)})^p
\end{eqnarray*}
modulo $p^{r-j+1}\Ld$. Set
\begin{eqnarray*}
G_{((p^{r-j-1}n)^{p^{i-r+j+1}})}^{(r)}&=&\frac1p(G_{((p^{r-j}n)^{p^{i-r+j}})}^{(r)}-(G_{((p^{r-j-1}n)^{p^{i-r+j}})}^{(r)})^p)\\
    &=&\hG_{((p^{r-j-1}n)^{p^{i-r+j}})}+p^{r-j}\ph_{i,j+1}(g).
\end{eqnarray*}
This completes the construction of $G_{p^rn}^{(r)}$ satisfying
properties (P3) and (P4).

\subsection{The $p^p$-form on $\Ld$}
Here we will extend the result one step further. There are two
reasons for considering this case. First, the nontrivial invariant
factor of the Cartan matrix of type $D_{2l+1}$ is 4. Therefore, by
Corollary \ref{cor:invariant factors}, we need to know the
invariant factors of the $2^2$-form for this case. Also, the
approach to this case lends insight into why the general proof is
so difficult.

We begin the construction the way we did in the previous section.
As before, we deduce that for $i<p$
\begin{eqnarray*}
G_{(p^in)}^{(p)}&=&\hG_{(p^in)}+p^{(p-1)(p-i)}(\hG_{(p^{i-1}n)})^p\\
        &=&p\hG_{((p^{i-1}n)^p)}+(1+p^{(p-1)(p-i)})(G_{(p^{i-1}n)}^{(p)}-p^{(p-1)(p-i+1)}(\hG_{(p^{i-2}n)})^p)^p.
\end{eqnarray*}
Note that when $i<p-1$, $(p-1)(p-i)>p-1>i$, hence
\[
G_{(p^in)}^{(p)}\equiv
p\hG_{((p^{i-1}n)^p)}+(G_{(p^{i-1}n)}^{(p)})^p
\]
modulo $p^{i+1}\Ld$. Therefore, in this case we can construct
$G_{(p^in)}^{(p)}$ satisfying property (P3) exactly as we did in
the previous section.

When $i=p-1$, we make the following adjustment. We have that
\begin{eqnarray*}
G_{(p^{p-1}n)}^{(p)}&=&\hG_{(p^{p-1}n)}+p^{p-1}(\hG_{(p^{p-2}n)})^p\\
        &=&p\hG_{((p^{p-2}n)^p)}+(1+p^{p-1})(G_{(p^{p-2}n)}^{(p)}-p^{2(p-1)}(\hG_{(p^{p-3}n)})^p)^p\\
        &\equiv&p\hG_{((p^{p-2}n)^p)}+(1+p^{p-1})(G_{(p^{p-2}n)}^{(p)})^p
\end{eqnarray*}
modulo $p^{p}\Ld$. Therefore, set
\begin{eqnarray*}
G_{((p^{p-2}n)^p)}^{(p)}&=&\frac1p[G_{(p^{p-1}n)}^{(p)}-(G_{(p^{p-2}n)}^{(p)})^p]\\
    &=&\hG_{((p^{p-2}n)^p)}+p^{p-2}(G_{(p^{p-2}n)}^{(p)})^p+p^{p-1}\ph_{p-1,1}(g).
\end{eqnarray*}

Next, assume by induction that $1\leq j<p-1$ and that we have
constructed
\[
G_{((p^{p-j-1})^{p^j})}^{(p)}=\hG_{((p^{p-j-1}n)^{p^j})}+p^{p-j-1}(G_{(p^{p-2}n)}^{(p)})^p+p^{p-j}\ph_{p-1,j}(g).
\]

By induction on $r$ ($=p$)
\[
\hG_{((p^{p-j-1}n)^{p^j})}=p\hG_{((p^{p-j-2}n)^{p^{j+1}})}+(\hG_{((p^{p-j-2}n)^{p^j})})^p,
\]
and (by induction on $i=p-1$)
\[
G_{((p^{p-j-2}n)^{p^j})}^{(p)}=\hG_{((p^{p-j-2}n)^{p^j})}+p^{p-j-1}\ph_{p-2,j}(g).
\]
Thus,
\begin{eqnarray*}
G_{((p^{p-j-1})^{p^j})}^{(p)}&=&p\hG_{((p^{p-j-2}n)^{p^{j+1}})}+p^{p-j-1}(G_{(p^{p-2}n)}^{(p)})^p\\
&&+(G_{((p^{p-j-2}n)^{p^j})}^{(p)}-p^{p-j-1}\ph_{p-2,j}(g))^p+p^{p-j}\ph_{p-1,j}(g)\\
    &\equiv&p\hG_{((p^{p-j-2}n)^{p^{j+1}})}+p^{p-j-1}(G_{(p^{p-2}n)}^{(p)})^p+(G_{((p^{p-j-2}n)^{p^j})}^{(p)})^p
\end{eqnarray*}
modulo $p^{p-j}\Ld$. Set
\begin{eqnarray*}
G_{((p^{p-j-2}n)^{p^{j+1}})}^{(p)}&=&\frac1p[G_{((p^{p-j-1})^{p^j})}^{(p)}-(G_{((p^{p-j-2}n)^{p^j})}^{(p)})^p]\\
    &=&\hG_{((p^{p-j-2}n)^{p^{j+1}})}+p^{p-j-2}(G_{(p^{p-2}n)}^{(p)})^p+p^{p-j-1}\ph_{p-1,j+1}(g).
\end{eqnarray*}
This completes the construction of $G_{(p^in)}^{p}$ for $i<p$
satisfying property (P3).

Assume that $i=p+k\geq p$. We will construct elements
$G_{((p^{p-j-1}n)^{p^{j+k+1}})}^{(p)}$ by induction on $k$ and
$j$. To this end, when $k=0$, we have
\begin{eqnarray*}
G_{(p^pn)}^{(p)}&=&p\hG_{(p^pn)}+(\hG_{(p^{p-1}n)})^p\\
        &=&p\hG_{(p^pn)}+(G_{(p^{p-1}n)}^{(p)}-p^{p-1}(\hG_{(p^{p-2}n)})^p)^p\\
        &=&p\hG_{(p^pn)}+(G_{(p^{p-1}n)}^{(p)}-p^{p-1}(G_{(p^{p-2}n)}^{(p)}-p^{2(p-1)}(\hG_{(p^{p-3}n)})^p)^p\\
        &\equiv& p\hG_{(p^pn)}+(G_{(p^{p-1}n)}^{(p)})^p
        -p^p\left(G_{(p^{p-1}n)}^{(p)}\right)^{p-1}\left(G_{(p^{p-2}n)}^{(p)}\right)^p
\end{eqnarray*}
modulo $p^{p+1}\Ld$. Set
\begin{eqnarray*}
G_{((p^{p-1}n)^p)}^{(p)}&=&\frac1p[G_{(p^pn)}^{(p)}-(G_{(p^{p-1}n)}^{(p)})^p]\\
    &=&\hG_{(p^pn)}-p^{p-1}\left(G_{(p^{p-1}n)}^{(p)}\right)^{p-1}\left(G_{(p^{p-2}n)}^{(p)}\right)^p
    +p^p\ph_{p,1}(g).
\end{eqnarray*}
Assume by induction that for $1\leq j<p-1$
\begin{eqnarray*}
G_{((p^{p-j}n)^{p^j})}^{(p)}&=&\hG_{((p^{p-j}n)^{p^j})}
-p^{p-j}\left[\sum_{l=0}^{j-1}\left(G_{((p^{p-l-1}n)^{p^l})}^{(p)}\right)^{p-1}\left(G_{(p^{p-2}n)}^{(p)}\right)^p\right]\\
&&+p^{p-j+1}\ph_{p,j}(g).
\end{eqnarray*}
By induction on $r$ ($=p$),
\[
\hG_{((p^{p-j}n)^{p^j})}=p\hG_{((p^{p-j-1}n)^{p^{j+1}})}+(\hG_{(p^{p-j-1}n)^{p^j}})^p,
\]
and, by induction on $i$ ($=p$),
\[
G_{((p^{p-j-1}n)^{p^j})}^{(p)}=\hG_{((p^{p-j-1}n)^{p^j})}+p^{p-j-1}(G_{(p^{p-2}n)}^{(p)})^p+p^{p-j}\ph_{p-1,j}(g).
\]
Hence
\begin{eqnarray*}
G_{((p^{p-j}n)^{p^j})}^{(p)}&=&p\hG_{((p^{p-j-1}n)^{p^{j+1}})}+(\hG_{((p^{p-j-1}n)^{p^j})})^p\\
&&-p^{p-j}\left[\sum_{l=0}^{j-1}\left(G_{((p^{p-l-1}n)^{p^l})}^{(p)}\right)^{p-1}\left(G_{(p^{p-2}n)}^{(p)}\right)^p\right]
+p^{p-j+1}\ph_{p,j}(g)\\
    &=&p\hG_{((p^{p-j-1}n)^{p^{j+1}})}\\
    &&+(G_{((p^{p-j-1}n)^{p^j})}^{(p)}-p^{p-j-1}(G_{(p^{p-2}n)}^{(p)})^p-p^{p-j}\ph_{p-1,j}(g))^p\\
&&-p^{p-j}\left[\sum_{l=0}^{j-1}\left(G_{((p^{p-l-1}n)^{p^l})}^{(p)}\right)^{p-1}\left(G_{(p^{p-2}n)}^{(p)}\right)^p\right]
+p^{p-j+1}\ph_{p,j}(g)\\
    &\equiv&p\hG_{((p^{p-j-1}n)^{p^{j+1}})}\\&&+(G_{((p^{p-j-1}n)^{p^j})}^{(p)})^p
    -p^{p-j}(G_{((p^{p-j-1}n)^{p^j})}^{(p)})^{p-1}(G_{(p^{p-2}n)}^{(p)})^p\\
    &&-p^{p-j}\left[\sum_{l=0}^{j-1}\left(G_{((p^{p-l-1}n)^{p^l})}^{(p)}\right)^{p-1}\left(G_{(p^{p-2}n)}^{(p)}\right)^p\right]\\
    &\equiv&p\hG_{((p^{p-j-1}n)^{p^{j+1}})}+(G_{((p^{p-j-1}n)^{p^j})}^{(p)})^p\\
    &&-p^{p-j}\left[\sum_{l=0}^{j}\left(G_{((p^{p-l-1}n)^{p^l})}^{(p)}\right)^{p-1}\left(G_{(p^{p-2}n)}^{(p)}\right)^p\right]
\end{eqnarray*}\
modulo $p^{p-j+1}\Ld$. Set
\begin{eqnarray*}
G_{((p^{p-j-1}n)^{p^{j+1}})}^{(p)}&=&\frac1p[G_{((p^{p-j}n)^{p^j})}^{(p)}-(G_{((p^{p-j-1}n)^{p^j})}^{(p)})^p]\\
    &=&\hG_{((p^{p-j-1}n)^{p^{j+1}})}+p^{p-j}\ph_{p,j+1}(g)\\
    &&-p^{p-j-1}\left[\sum_{l=0}^{j}\left(G_{((p^{p-l-1}n)^{p^l})}^{(p)}\right)^{p-1}
    \left(G_{(p^{p-2}n)}^{(p)}\right)^p\right].
\end{eqnarray*}
When $j=p-2$, we have
\begin{eqnarray*}
G_{((pn)^{p^{p-1}})}^{(p)}&=&\hG_{((pn)^{p^{p-1}})}
    -p\left[\sum_{l=0}^{p-2}\left(G_{((p^{p-l-1}n)^{p^l})}^{(p)}\right)^{p-1}\left(G_{(p^{p-2}n)}^{(p)}\right)^p\right]\\
    &&+p^2\ph_{p,p-1}(g).
\end{eqnarray*}
As before,
\[
\hG_{((pn)^{p^{p-1}})}=p\hG_{(n^{p^p})}+(\hG_{(n^{p^{p-1}})})^p,
\]
but this time
\begin{eqnarray*}
G_{(n^{p^{p-1}})}^{(p)}&=&\hG_{(n^{p^{p-1}})}+(G_{(p^{p-2}n)}^{(p)})^p+p\ph_{p-1,p-1}(g).
\end{eqnarray*}
Therefore,
\begin{eqnarray*}
G_{((pn)^{p^{p-1}})}^{(p)}&=&p\hG_{(n^{p^p})}+(\hG_{(n^{p^{p-1}})})^p\\
    &&-p\left[\sum_{l=0}^{p-2}\left(G_{((p^{p-l-1}n)^{p^l})}^{(p)}\right)^{p-1}\left(G_{(p^{p-2}n)}^{(p)}\right)^p\right]
    +p^2\ph_{p,p-1}(g)\\
    &=&p\hG_{(n^{p^p})}+(G_{(n^{p^{p-1}})}^{(p)}-(G_{(p^{p-2}n)}^{(p)})^p-p\ph_{p-1,p-1}(g))^p\\
    &&-p\left[\sum_{l=0}^{p-2}\left(G_{((p^{p-l-1}n)^{p^l})}^{(p)}\right)^{p-1}\left(G_{(p^{p-2}n)}^{(p)}\right)^p\right]
    +p^2\ph_{p,p-1}(g)\\
    &=&p\hG_{(n^{p^p})}+(G_{(n^{p^{p-1}})}^{(p)}-(G_{(p^{p-1}n)}^{(p)}-pG_{((p^{p-2}n)^p)}^{(p)})-p\ph_{p-1,p-1}(g))^p\\
    &&-p\left[\sum_{l=0}^{p-2}\left(G_{((p^{p-l-1}n)^{p^l})}^{(p)}\right)^{p-1}\left(G_{(p^{p-2}n)}^{(p)}\right)^p\right]
    +p^2\ph_{p,p-1}(g)\\
    &\equiv&p\hG_{(n^{p^p})}+(G_{(n^{p^{p-1}})}^{(p)}-G_{(p^{p-1}n)}^{(p)})^p\\
    &&-p\left[\sum_{l=0}^{p-2}\left(G_{((p^{p-l-1}n)^{p^l})}^{(p)}\right)^{p-1}\left(G_{(p^{p-2}n)}^{(p)}\right)^p\right]
\end{eqnarray*}
modulo $p^2\Ld$. Set
\begin{eqnarray*}
G_{(n^{p^p})}^{(p)}&=&\frac1p[G_{((pn)^{p^{p-1}})}^{(p)}-(G_{(n^{p^{p-1}})}^{(p)}-G_{(p^{p-1}n)}^{(p)})^p]\\
    &=&\hG_{(n^{p^p})}
    -\left[\sum_{l=0}^{p-2}\left(G_{((p^{p-l-1}n)^{p^l})}^{(p)}\right)^{p-1}\left(G_{(p^{p-2}n)}^{(p)}\right)^p\right]\\
    &&+p\ph_{p,p}(g).
\end{eqnarray*}


We are now ready to complete the result. To this end, assume that
$k> 0$ and that
\begin{small}
\begin{enumerate}
\item[(I1)]
\begin{eqnarray*}
G_{((p^{p-1}n)^{p^{k}})}^{(p)}&=&\hG_{(p^{p+k-1}n)}+p^p\ph_{p+k-1,k}(g)\\
    &&+(-1)^{k}
    p^{p-1}\prod_{m=0}^{k-1}\left(G_{((p^{p-1}n)^{p^m})}^{(p)}\right)^{p-1}\left(G_{(p^{p-2}n)}^{(p)}\right)^p.
\end{eqnarray*}
\item[(I2)] For $j\geq 1$,
\begin{eqnarray*}
G_{((p^{p-j-1}n)^{p^{j+k}})}^{(p)}&=&\hG_{((p^{p-j-1}n)^{p^{j+k}})}+p^{p-j}\ph_{p+k-1,j+k}(g)\\
    &&\hspace{-.9in}+(-1)^{k}(p^{p-j-1})\sum_{0\leq l_0\leq\cdots\leq l_{k-1}\leq j}\left[\prod_{m=0}^{k-1}
        \left(G_{((p^{p-l_m-1}n)^{p^{l_m+m}})}^{(p)}\right)^{p-1}\left(G_{(p^{p-2}n)}^{(p)}\right)^p\right].
\end{eqnarray*}
\end{enumerate}
\end{small}
Notice that when $k=1$ this agrees with the previous case.

Now, observe that
\begin{eqnarray*}
G_{(p^{p+k}n)}^{(p)}&=&g_{p^{p+k}n}^{(p,p)}\\
            &=&pg_{p^{p+k}n}^{(p-1,p)}+(g_{p^{p+k-1}n}^{(p-1,p)})^p\\
            &=&p\hG_{(p^{p+k}n)}+(\hG_{(p^{p+k-1}n)})^p\\
            &=&p\hG_{(p^{p+k}n)}+
            \Bigg(G_{((p^{p-1}n)^{p^{k}})}^{(p)}-(-1)^{k}p^{p-1}\\
            &&\left.\times\prod_{m=0}^{k-1}\left(G_{((p^{p-1}n)^{p^m})}^{(p)}\right)^{p-1}\left(G_{(p^{p-2}n)}^{(p)}\right)^p
    +p^p\ph_{p+k-1,k}(g)\right)^p\\
            &\equiv&p\hG_{(p^{p+k}n)}+(G_{((p^{p-1}n)^{p^{k}})}^{(p)})^p\\
    &&+(-1)^{k+1}p^p\left(G_{((p^{p-1}n)^{p^{k}})}^{(p)}\right)^{p-1}
    \prod_{m=0}^{k-1}\left(G_{((p^{p-1}n)^{p^m})}^{(p)}\right)^{p-1}\left(G_{(p^{p-2}n)}^{(p)}\right)^p\\
    &\equiv&p\hG_{(p^{p+k}n)}+(G_{((p^{p-1}n)^{p^{k}})}^{(p)})^p\\
    &&+(-1)^{k+1}p^p\prod_{m=0}^{k}\left(G_{((p^{p-1}n)^{p^m})}^{(p)}\right)^{p-1}\left(G_{(p^{p-2}n)}^{(p)}\right)^p\\
\end{eqnarray*}
modulo $p^{p+1}\Ld$. Set
\begin{eqnarray*}
G^{(p)}_{((p^{p-1}n)^{p^{k+1}})}&=&\frac1p[G_{(p^{p+k}n)}^{(p)}-(G_{((p^{p-1}n)^{p^{k}})}^{(p)})^p]\\
    &=&\hG_{(p^{p+k}n)}+
    (-1)^{k+1}p^{p-1}\prod_{m=0}^{k}\left(G_{((p^{p-1}n)^{p^m})}^{(p)}\right)^{p-1}\left(G_{(p^{p-2}n)}^{(p)}\right)^p\\
    &&+p^p\ph_{p+k,k+1}(g).
\end{eqnarray*}
Next assume that $1\leq j<p-1$ and
\begin{small}
\begin{eqnarray*}
G^{(p)}_{((p^{p-j}n)^{p^{j+k}})}&=&\hG_{((p^{p-j}n)^{p^{j+k}})}+p^{p-j+1}\ph_{p+k,j+k}(g)\\
    &&\hspace{-.75in}+(-1)^{k+1}(p^{p-j})\sum_{0\leq l_0\leq\cdots\leq l_{k}\leq j-1}\left[\prod_{m=0}^{k}
        \left(G_{((p^{p-l_m-1}n)^{p^{l_m+m}})}^{(p)}\right)^{p-1}\left(G_{(p^{p-2}n)}^{(p)}\right)^p\right].
\end{eqnarray*}
\end{small}
By induction on $r$,
\[
\hG_{((p^{p-j}n)^{p^{j+k}})}=p\hG_{((p^{p-j-1}n)^{p^{j+k+1}})}+(\hG_{((p^{p-j-1}n)^{p^{j+k}})})^p,
\]
and, by induction on $k$,
\begin{tiny}
\begin{eqnarray*}
G_{((p^{p-j-1}n)^{p^{j+k}})}^{(p)}&=&\hG_{((p^{p-j-1}n)^{p^{j+k}})}+p^{p-j}\ph_{p+k-1,j+k}(g)\\
    &&\hspace{-.75in}+(-1)^{k}(p^{p-j-1})\sum_{0\leq l_0\leq\cdots\leq l_{k-1}\leq j}\left[\prod_{m=0}^{k-1}
        \left(G_{((p^{p-l_m-1}n)^{p^{l_m+m}})}^{(p)}\right)^{p-1}\left(G_{(p^{p-2}n)}^{(p)}\right)^p\right].
\end{eqnarray*}
\end{tiny}
Hence, modulo $p^{p-j+1}\Ld$,
\begin{tiny}
\begin{eqnarray*}
G^{(p)}_{((p^{p-j}n)^{p^{j+k}})}&=&p\hG_{((p^{p-j-1}n)^{p^{j+k+1}})}+(\hG_{((p^{p-j-1}n)^{p^{j+k}})})^p+p^{p-j+1}\ph_{p+k,j+k}(g)\\
    &&\hspace{-.75in}+(-1)^{k+1}(p^{p-j})\sum_{0\leq l_0\leq\cdots\leq l_{k}\leq j}\left[\prod_{m=0}^{k}
        \left(G_{((p^{p-l_m-1}n)^{p^{l_m+m}})}^{(p)}\right)^{p-1}\left(G_{(p^{p-2}n)}^{(p)}\right)^p\right]\\
    &\equiv&p\hG_{((p^{p-j-1}n)^{p^{j+k+1}})}+(G_{((p^{p-j-1}n)^{p^{j+k}})}^{(p)})^p\\
    &&\hspace{-1in}-(-1)^{k}p^{p-j}(G_{((p^{p-j-1}n)^{p^{j+k}})}^{(p)})^{p-1}
    \sum_{0\leq l_0\leq\cdots\leq l_{k-1}\leq j}\left[\prod_{m=0}^{k-1}
        \left(G_{((p^{p-l_m-1}n)^{p^{l_m+m}})}^{(p)}\right)^{p-1}\left(G_{(p^{p-2}n)}^{(p)}\right)^p\right]\\
    &&\hspace{-.75in}+(-1)^{k+1}(p^{p-j})
        \sum_{0\leq l_0\leq\cdots\leq l_{k}\leq j-1}\left[\prod_{m=0}^{k}
        \left(G_{((p^{p-l_m-1}n)^{p^{l_m+m}})}^{(p)}\right)^{p-1}\left(G_{(p^{p-2}n)}^{(p)}\right)^p\right]\\
    &\equiv&p\hG_{((p^{p-j-1}n)^{p^{j+k+1}})}+(G_{((p^{p-j-1}n)^{p^{j+k}})}^{(p)})^p\\
    &&\hspace{-.5in}+(-1)^{k+1}(p^{p-j})
        \sum_{0\leq l_0\leq\cdots\leq l_{k}\leq j}\left[\prod_{m=0}^{k}
        \left(G_{((p^{p-l_m-1}n)^{p^{l_m+m}})}^{(p)}\right)^{p-1}\left(G_{(p^{p-2}n)}^{(p)}\right)^p\right].
\end{eqnarray*}
\end{tiny}
Set
\begin{tiny}
\begin{eqnarray*}
G_{((p^{p-j-1}n)^{p^{j+k+1}})}&=&\frac1p[G_{((p^{p-j}n)^{p^{j+k}})}-(G_{((p^{p-j-1}n)^{p^{j+k}})}^{(p)})^p]\\
    &=&\hG_{((p^{p-j-1}n)^{p^{j+k+1}})}+p^{p-j}\ph_{p+k,j+k+1}(g)\\
    &&\hspace{-.75in}+(-1)^{k+1}(p^{p-j-1})
    \sum_{0\leq l_0\leq\cdots\leq l_{k}\leq j}\left[\prod_{m=0}^{k}
        \left(G_{((p^{p-l_m-1}n)^{p^{l_m+m}})}^{(p)}\right)^{p-1}\left(G_{(p^{p-2}n)}^{(p)}\right)^p\right].
\end{eqnarray*}
\end{tiny}
We have verified the inductive assumptions above!

When $j=p-2$, we have
\begin{eqnarray*}
G_{((pn)^{p^{p+k-1}})}^{(p)}&=&\hG_{((pn)^{p^{p+k-1}})}+p^{2}\ph_{p+k,p+k-1}(g)\\
    &&\hspace{-.75in}+(-1)^{k+1} p\sum_{0\leq l_0\leq\cdots\leq l_{k}\leq p-2}\left[\prod_{m=0}^{k}
        \left(G_{((p^{p-l_m-1}n)^{p^{l_m+m}})}^{(p)}\right)^{p-1}\left(G_{(p^{p-2}n)}^{(p)}\right)^p\right].
\end{eqnarray*}
As usual
\[
\hG_{((pn)^{p^{p+k-1}})}=p\hG_{(n^{p^{p+k}})}+(\hG_{(n^{p^{p+k-1}})})^p
\]
But this time
\begin{eqnarray*}
G_{(n^{p^{p+k-1}})}^{(p)}&=&\hG_{(n^{p^{p+k-1}})}+p\ph_{p+k-1,p+k-1}(g)\\
    &&\hspace{-.75in}+(-1)^{k}\sum_{0\leq l_0\leq\cdots\leq l_{k-1}\leq p-1}\left[\prod_{m=0}^{k-1}
        \left(G_{((p^{p-l_m-1}n)^{p^{l_m+m}})}^{(p)}\right)^{p-1}\left(G_{(p^{p-2}n)}^{(p)}\right)^p\right].
\end{eqnarray*}
Hence
\begin{tiny}
\begin{eqnarray*}
\hspace{-.25in}G_{((pn)^{p^{p+k-1}})}^{(p)}&\equiv&p\hG_{(n^{p^{p+k}})}+p\sum_{0\leq
l_0\leq\cdots\leq l_{k}\leq p-2}\left[\prod_{m=0}^{k}
        \left(G_{((p^{p-l_m-1}n)^{p^{l_m+m}})}^{(p)}\right)^{p-1}\left(G_{(p^{p-2}n)}^{(p)}\right)^p\right]\\
    &&\hspace{-1in}+
    \left(G_{((pn)^{p^{p+k-1}})}^{(p)}-(-1)^{k}\sum_{0\leq l_0\leq\cdots\leq l_{k-1}\leq p-1}\left[\prod_{m=0}^{k-1}
        \left(G_{((p^{p-l_m-1}n)^{p^{l_m+m}})}^{(p)}\right)^{p-1}\left(G_{(p^{p-2}n)}^{(p)}\right)^p\right]\right)^p\\
    &\equiv&p\hG_{(n^{p^{p+k}})}+p\sum_{0\leq l_0\leq\cdots\leq l_{k}\leq p-2}\left[\prod_{m=0}^{k}
        \left(G_{((p^{p-l_m-1}n)^{p^{l_m+m}})}^{(p)}\right)^{p-1}\left(G_{(p^{p-2}n)}^{(p)}\right)^p\right]\\
    &&\hspace{-1in}+
    \left(G_{((pn)^{p^{p+k-1}})}^{(p)}-(-1)^{k}\sum_{0\leq l_0\leq\cdots\leq l_{k-1}\leq p-1}\left[\prod_{m=0}^{k-1}
        \left(G_{((p^{p-l_m-1}n)^{p^{l_m+m}})}^{(p)}\right)^{p-1}\left(G_{(p^{p-1}n)}^{(p)}\right)\right]\right)^p.
\end{eqnarray*}
\end{tiny}
modulo $p^2\Ld$. (Note that, in the last line of the computation
above, we have used $(G_{(p^{p-2}n)}^{(p)})^p\equiv
G_{(p^{p-1}n)}^{(p)}$ modulo $p\Ld$ and the fact that
$(x+py)^p\equiv x$ mod $p^2$.)

\bigskip

Set
\begin{tiny}
\begin{eqnarray*}
G_{(n^{p^{p+k}})}^{(p)}&=&
    \frac1p\Bigg[ G_{((pn)^{p^{p+k-1}})}^{(p)}-
    \Bigg( G_{((pn)^{p^{p+k-1}})}^{(p)}+(-1)^{k+1}\\
    &&\times\sum_{0\leq l_0\leq\cdots\leq l_{k-1}\leq p-1}\prod_{m=0}^{k-1}
        \left(G_{((p^{p-l_m-1}n)^{p^{l_m+m}})}^{(p)}\right)^{p-1}\left(G_{(p^{p-1}n)}^{(p)}\right)\Bigg)^p\Bigg]
        \\
    &=&\hG_{(n^{p^{p+k}})}+\sum_{0\leq l_0\leq\cdots\leq l_{k}\leq p-2}\left[\prod_{m=0}^{k}
        \left(G_{((p^{p-l_m-1}n)^{p^{l_m+m}})}^{(p)}\right)^{p-1}\left(G_{(p^{p-2}n)}^{(p)}\right)^p\right]\\
        &&+p\ph_{p+k,p+k}(g).
\end{eqnarray*}
\end{tiny}
This completes the construction of $(G_\ld^{(p)})$.


\subsection{Statement of the Main Result}

\begin{thm}\label{thm:r<p} Let $r\leq p$. Then, there exists a bases
$(G_{\ld}^{(r)})$ and $(N_\ld^{(r)})$ for $\Ld$ such that
\[
\Pi(x,y)^{p^r}=
    \sum_\ld D_r(\ld)G_\ld^{(r)}(x)N_\ld^{(r)}(y).
\]
\end{thm}

\end{document}